\documentclass[conference]{IEEEtran}

\ifCLASSINFOpdf
   \usepackage[pdftex]{graphicx}
   \DeclareGraphicsExtensions{.pdf,.jpeg,.png}
\else
\fi
\usepackage{amsthm}
\newtheorem*{remark}{Remark}
\newtheorem{definition}{Definition}[section]

\usepackage{xcolor}
\usepackage{soul}
\usepackage{amsmath}
\usepackage{amsfonts}
\usepackage{booktabs}
\usepackage{url}

\usepackage[]{fancyhdr} % 
\newcommand{\changefont}{\fontsize{9}{9}\selectfont}
\ifCLASSOPTIONcompsoc
  \usepackage[caption=false,font=normalsize,labelfont=sf,textfont=sf]{subfig}
\else
  \usepackage[caption=false,font=footnotesize]{subfig}
\fi
\IEEEoverridecommandlockouts
\begin{document}
\title{
Exploring reactive power limits on wind\\ farm collector networks with \\convex inner approximations
}

\author{\IEEEauthorblockN{Nawaf Nazir}
\IEEEauthorblockA{\textit{Energy and Environment Directorate}\\Pacific Northwest National Laboratory\\
Richland, USA\\
nawaf.nazir@pnnl.gov}
\and
\IEEEauthorblockN{Ian A. Hiskens}
\IEEEauthorblockA{\textit{Department of EECS}\\University of Michigan\\
Ann Arbor, USA\\
hiskens@umich.edu}
\and
\IEEEauthorblockN{Mads R. Almassalkhi}
\IEEEauthorblockA{\textit{Department of EBE}\\University of Vermont\\
Burlington, USA\\
malmassa@uvm.edu}
\thanks{N. Nazir and M. Almassalkhi were supported by the U.S.
Department of Energy’s Advanced Research Projects Agency-Energy Award DE-AR0000694 and the National Science Foundation (NSF) Award ECCS-2047306.}}


%Telephone: (800) 555--1212, Fax: (888) 555--1212}

% <-this % stops a space


%\lhead{11TH BULK POWER SYSTEMS DYNAMICS AND CONTROL SYMPOSIUM, JULY 25-30, 2022, BANFF, CANADA}
%\rhead{1}

%\fontfamily{phv}\fontseries{b}\fontsize{9}{11}\selectfont

\maketitle
\thispagestyle{fancy}
\pagestyle{fancy}

%\thispagestyle{fancy}
%\pagestyle{fancy}

%%%%%%%%%%%%%%%%%%%%%%%%%%%%%%%%%%%%%%%%%%%%%%%%%%%%%%%%%%%%%%%%%%%%%%%%%%%%%%%%
\begin{abstract}

A wind farm can provide reactive power at sub-transmission and transmission buses in order to support and improve voltage profiles. It is common for the reactive power capability of a wind farm to be evaluated as the sum of the individual turbine ratings. However, such an assessment does not take into account losses over the collector network, nor the voltage constraints imposed by the turbines and network. In contrast, the paper presents a method for determining the range of reactive power support that each turbine can provide whilst guaranteeing satisfaction of voltage constraints. This is achieved by constructing convex inner approximations of the non-convex set of admissible reactive power injections. We present theoretical analysis that supports the constraint satisfaction guarantees. An example illustrates the effectiveness of the algorithm and provides a comparison with a fully decentralized approach to controlling wind farm reactive power. Such approaches have the potential to improve the design and operation of wind farm collector networks, reducing the need for additional costly reactive power resources.

\end{abstract}

%%%%%%%%%%%%%%%%%%%%%%%%%%%%%%%%%%%%%%%%%%%%%%%%%%%%%%%%%%%%%%%%%%%%%%%%%%%%%%%%
\section{INTRODUCTION}
Reactive power support from wind farms can play on important role in maintaining power system reliability. Specifically, the reactive power capability of type~3 and type~4 wind turbines can be used to regulate the grid voltage at the point of common coupling (PCC)~\cite{slootweg2005wind}. Hence, it is important for wind farm operators to characterize and control their reactive power capability so that this resource is available to the transmission system operator (TSO). Early work
characterizing the reactive power capability did not account for the collector network that interconnects the wind turbines~\cite{Camm_2009}. Consequently, the impact of voltage limits could not be assessed~\cite{opila2010wind}. This issue was partially addressed by the decentralized control scheme proposed in~\cite{hiskens2013strategies,Martins_2015}, which controls wind turbine reactive power to regulate the PCC voltage but does not offer {\em a priori} assessment of the available reactive power capability. More recent work provides voltage support from wind farms by using a sensitivity-based approach to rank reactive power loading for wind farms and their turbines~\cite{silva2019loading}.

A alternative approach to account for the wind farm network while dispatching turbine reactive power, is to explicitly consider the wind farm's radial (balanced) network within an optimal power flow (OPF) setting. However, that requires solving a non-convex OPF problem, which is NP-hard~\cite{molzahn2017computing}. The technical challenges associated with the non-convex formulation could be overcome by considering either linear approximations or convex relaxations~\cite{molzahn2019survey}. For example, traditional methods for solving the OPF problem in (radial) networks include the \textit{LinDist} model, which neglects the losses in the network to arrive at a simplified linear model. Much of the previous work on wind farm optimization utilizes the \textit{LinDist} model as it offers computational benefits. However, ignoring line losses (both reactive and active) can lead to unmodeled voltage violations under certain operating conditions~\cite{nazir2019convex}. Since convex relaxations of the OPF problem can provide solutions with zero duality~\cite{gan2014exact}, they have become popular proxies for the underlying network physics. However, in the case of a wind farm providing a desired value of reactive power at the PCC, convex relaxations can engender optimal solutions with so-called fictitious losses whose realized dispatch can cause voltages to exceed their limits~\cite{brahma2020vb}.

This paper overcomes previous shortcomings by employing convex inner approximations (or convex restrictions) to determine a wind farm's realizable reactive power capacity, and devise a feedback control scheme for regulating the PCC voltage. The control strategy disaggregates the time-varying reactive power reference among the wind turbines in a manner that guarantees network conditions always remain within limits. Unlike convex relaxations (i.e., outer approximations) and linearized approximations, convex inner approximations (CIAs) ensure that feasible solutions are also physically realizable. Of course, inner approximation may beget conservative solutions that can reduce performance.

Previously, CIAs have been employed in the optimization of dispatching (discrete) mechanical grid assets~\cite{nazir2019voltage} and (continuous) distributed energy resources (DERs)~\cite{nazir2019grid}. In this work, we adapt CIAs to determine practical reactive power bounds for each turbine (i.e., at each node) for a given wind power scenario. Within these nodal reactive power bounds, we can guarantee that any combination of turbine reactive power dispatch will ensure that network voltages are within their limits. Based on these nodal bounds, a real-time disaggregation control loop is formulated that can dispatch turbines and deliver desired reactive power to support grid operations.

The paper is organised as follows: Section~\ref{sec:model} develops the mathematical model of a wind farm network and illustrates the concept of nodal reactive capacities on a simple 3-node wind farm. Section~\ref{sec:ciaOpti} develops the convex inner approximation for  the non-convex optimization problem that defines the wind farm's reactive power capability. Section~\ref{sec:rtAVRdisagg} develops a real-time control algorithm that provides grid voltage support while ensuring satisfaction of wind farm network voltage constraints. Section~\ref{sec:concl} concludes the paper and highlights future research directions.

\section{Mathematical modeling and nodal reactive capacity}\label{sec:model}

\subsection{Wind farm model}

\begin{figure}[t]
\centering
\includegraphics[width=0.38\textwidth]{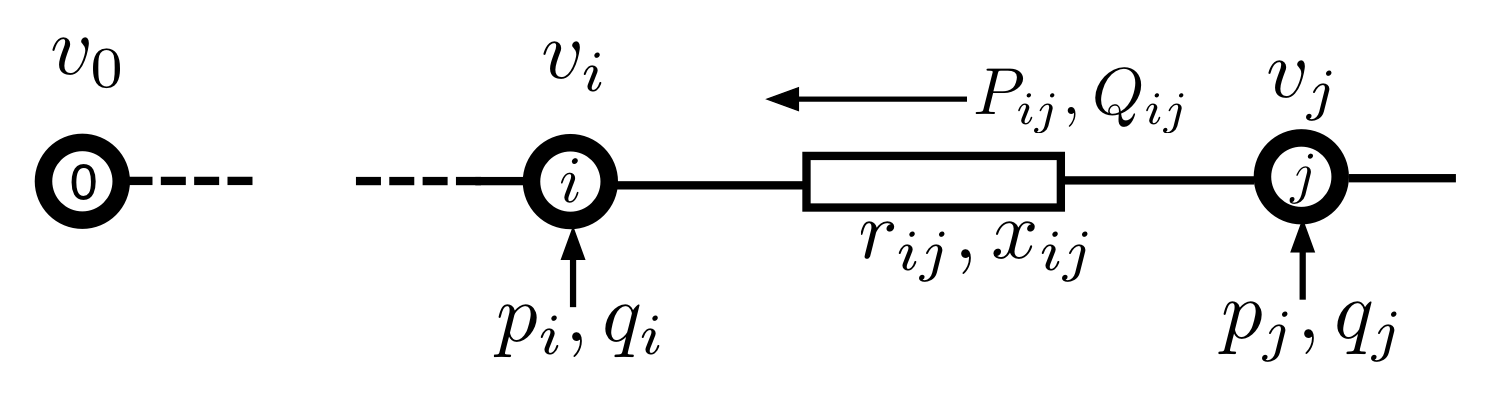}
\caption{\label{fig:radial_network} Nomenclature for a radial wind farm network~\cite{heidari2017non}.}
\end{figure}

In this section, we present the model of a wind farm, where a balanced, radial network often couples the turbines to the PCC, as shown in Fig.~\ref{fig:radial_network}. Thus, we can use the nonlinear \textit{DistFlow} formulation to model the wind farm network. Consider an undirected graph $\mathcal{G}=\{\mathcal{N}\cup\{0\}, \mathcal{L}\}$ consisting of a set of $N+1$ nodes with $\mathcal{N}:=\{1, \hdots, N\}$ and a set of $N$ branches $\mathcal{L}:=\{1, \hdots, N\} \subseteq \mathcal{N}\times \mathcal{N}$, such that $(i,j)\in \mathcal{L}$, if nodes $i,j$ are connected. %, where $L=N$ for a radial network. 
Node $0$ is assumed to be the head node (i.e., PCC) with a fixed voltage $V_0$.
Let $B \in \mathbb{R}^{(N+1) \times N}$ be the \textit{incidence matrix} of $\mathcal{G}$ relating the branches in $\mathcal{L}$ to the nodes in $\mathcal{N}\cup \{0\}$, such that the $(i,k)$-th entry of $B$ is $1$ if the $i$-th node is connected to the $k$-th branch and, otherwise, $0$. Without loss of generality, $B$ can be organized to form an upper-triangular matrix. If $V_i$ and $V_j$ are the voltage phasors at nodes $i$ and $j$ and $I_{ij}$ is the current phasor in branch $(i,j)\in \mathcal{L}$, then define $v_i:=|V_i|^2$, $v_j:=|V_j|^2$ and $l_{ij}:=|I_{ij}|^2$. Let $P_{ij}$ ($Q_{ij}$) be the active (reactive) power flow from node $j$ to $i$, let $p_j$ ($q_j$) be the active (reactive) power generations into node $j$, and let $r_{ij}$ ($x_{ij}$) be the resistance (reactance) of branch $(i,j)\in \mathcal{L}$, which means that the branch impedance is given by $z_{ij}:=r_{ij}+jx_{ij}$. Then, for a radial wind farm, the relation between node voltages and power flows is given by the \textit{DistFlow} equations $\forall (i,j)\in \mathcal{L}$:
\begin{subequations}\label{eq:dist_flow}
\begin{align}
v_j=&v_i+2r_{ij}P_{ij}+2x_{ij}Q_{ij}-|z_{ij}|^2l_{ij} \label{eq:volt_rel}\\
P_{ij}=&p_j+\sum_{h:h\rightarrow j}(P_{jh}-r_{jh}l_{jh}) \label{eq:real_power_rel}\\
Q_{ij}=&q_j+\sum_{h:h\rightarrow j}(Q_{jh}-x_{jh}l_{jh}) \label{eq:reac_power_rel}\\
l_{ij}(P_{ij},Q_{ij},v_j)  =& \frac{P_{ij}^2+Q_{ij}^2}{v_j}, \label{eq:curr_rel}
\end{align}
\end{subequations}
%where nodal power injections are $p_j:=p_{\text{g},j}-P_{\text{L},j}$ and $q_j:=q_{\text{g},j}-Q_{\text{L},j}$ with $p_{\text{g},j}$ ($q_{\text{g},j}$) as the controllable active (reactive) injections and $P_{\text{L},j}$ ($Q_{\text{L},j}$) is the uncontrollable active (reactive) demand data. 

The goal of this work is to maximize the range of reactive power output from the wind farm, i.e., $Q_{10}$, such that all  voltages $v_j$ and currents $l_{ij}$ are within their respective limits (i.e., $v_j \in [\underline{v}_j,\overline{v}_j] \,\, \forall j\in \mathcal{N}$ and $l_{ij} \in [\underline{l}_{ij},\overline{l}_{ij}] \,\, \forall (i,j)\in \mathcal{L}$).
%The goal is then to find the range of controllable variables (e.g., $[\underline{p},\overline{p}]$) within which all constraint limits are satisfied, i.e., the hosting capacity. 
However, finding such a range is challenging due to the non-linear nature of~\eqref{eq:curr_rel}. For clarity, we provide definitions of the following key terms used in the manuscript.
\begin{definition}[AC Admissibility]
A solution of a convex OPF problem is AC admissible, if the solution applied to the original, non-convex AC OPF, which uses~\eqref{eq:dist_flow}, is feasible. %, i.e., $v_j\in [\underline{V},\overline{V}]\,\, \forall j\in\mathcal{N}$ and $l_{ij}\in [\underline{l},\overline{l}] \,\,\forall (i,j) \in \mathcal{L}$.% are within their limits, i.e., $[\underline{V},\overline{V}]$ and $[\underline{l},\overline{l}]$
\end{definition}
\begin{definition}[Nodal reactive capacity]
Nodal reactive capacity is the range of AC admissible reactive power dispatch $\Delta q_{j} := [q_j^-, q_j^+]\, \forall j\in \mathcal{N}$ with lower and upper bounds $q_{j}^-\le 0$ and $q_{j}^+\ge 0$, respectively. That is, for all nodes $j$, any dispatch $q_{j} \in \Delta q_{j}$ is AC admissible.
\end{definition}

Next, we consider the nodal reactive capacity of a simple 3-node wind farm network in Fig.~\ref{fig:3_node_model} to motivate the approach.

\subsection{Motivating example on nodal reactive capacity}\label{sec:3-node_set}

\begin{figure}[t]
\centering
\includegraphics[width=0.22\textwidth]{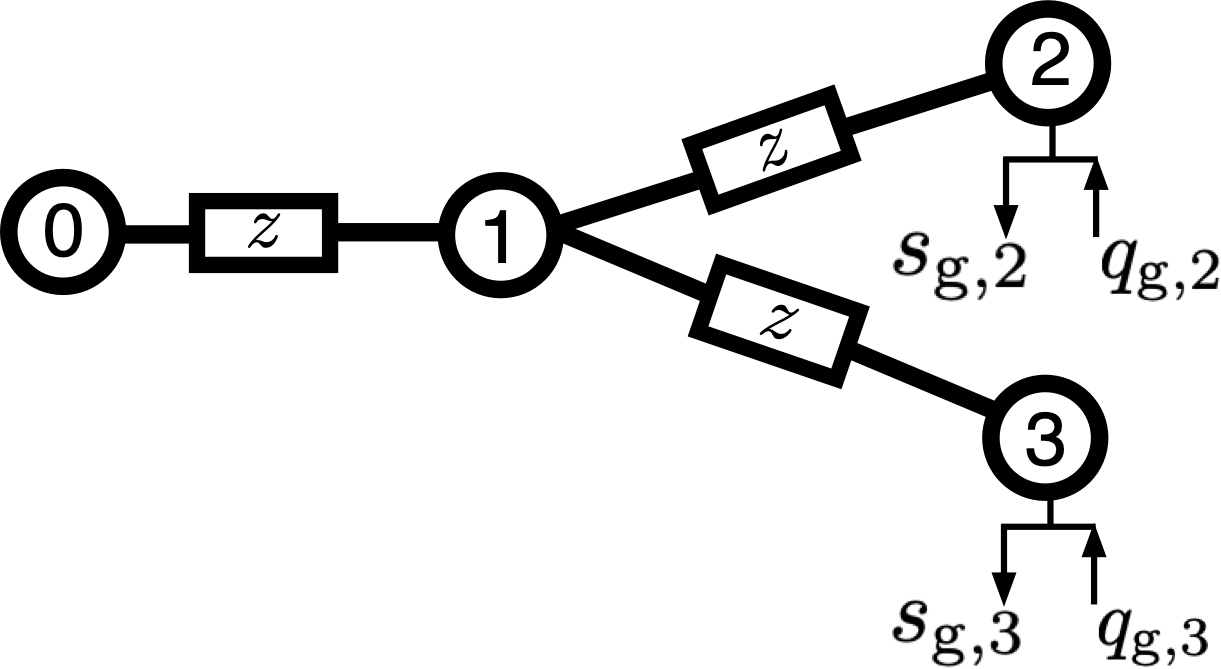}
\caption{\label{fig:3_node_model} A simple 3-node wind farm network.}
\end{figure}

\begin{figure}[t]
\centering
\includegraphics[width=0.36\textwidth]{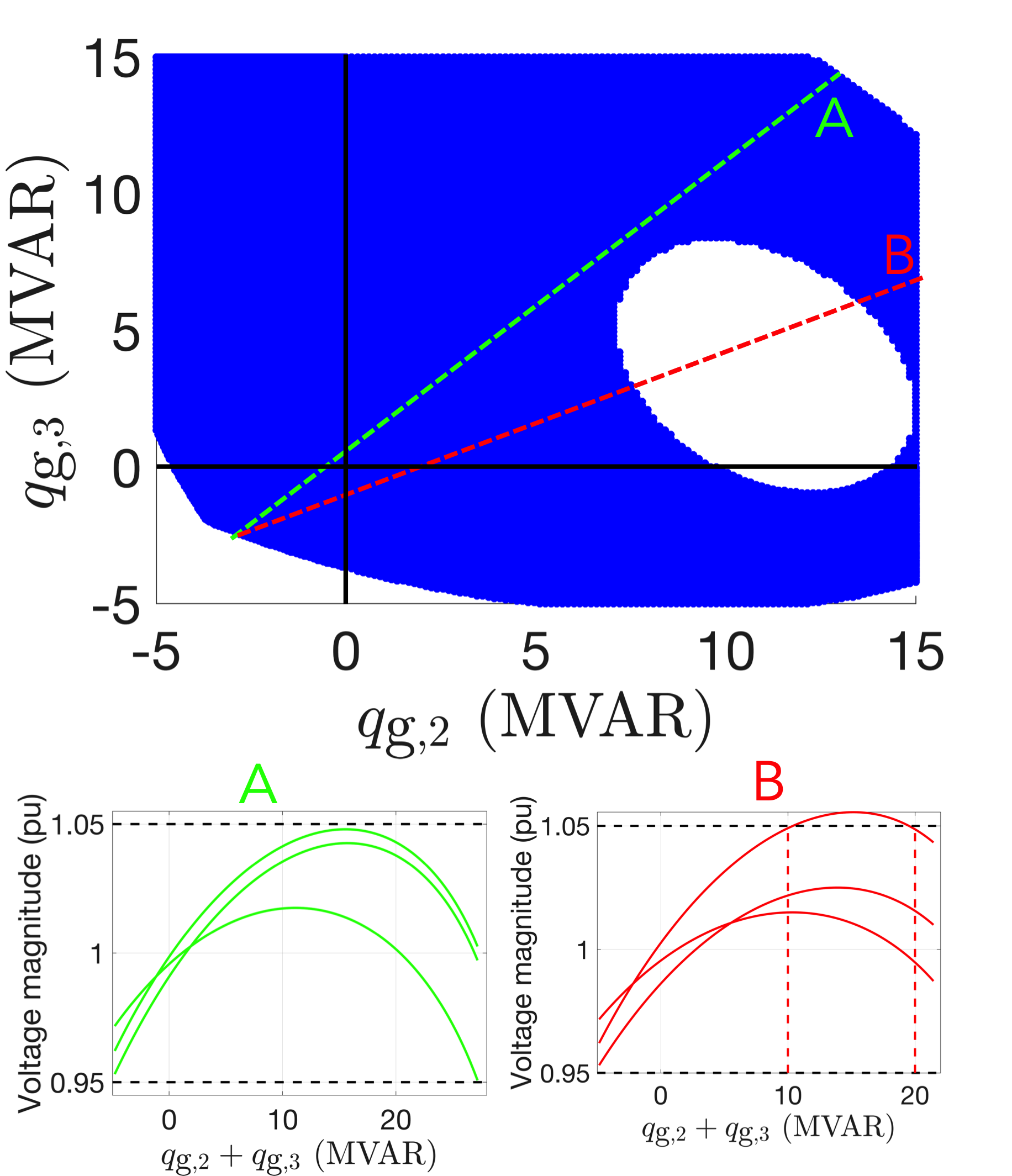}
\caption{\label{fig:3_node_sweep} Analysis of 3-node example. (Top) The set of admissible injections is non-convex. (Bottom) Voltage profiles along admissible (A) and inadmissible (B) trajectories
}
\end{figure}

Fig.~\ref{fig:3_node_model} represents a simple, balanced wind farm network with two turbines at nodes $2$ and $3$, and $V_0=1$~pu. Each (positive sequence) branch of the network has impedance $z=0.228+0.092j$~pu. Nodes~$2$ and $3$ have generation $s_{\text{g},2}=0.005-0.02j$ pu and $s_{\text{g},3}=0.01-0.015j$ pu, respectively.
%The base rating of the system is $100$MVA.
%Flexible resources $q_{\text{g},2}$ and $q_{\text{g},3}$ are assumed to be located at nodes $2$ and $3$. 
Only the reactive power injections at nodes $2$ and $3$ (labelled $q_{\text{g},2}$ and $q_{\text{g},3}$) are assumed to be controllable. Based on the AC power flow solutions obtained with Matpower~\cite{zimmerman2011matpower}, by varying $q_{\text{g},2}$ and $q_{\text{g},3}$, Fig.~\ref{fig:3_node_sweep} shows the feasible set of the AC OPF for the 3-node system. The figure shows that the admissible set is non-convex and contains a ``hole'' due to a voltage constraint. Hence, it is important when dispatching $q_{\text{g},2}$ and $q_{\text{g},3}$ to choose a trajectory that ensures AC admissibility. Specifically, Fig.~\ref{fig:3_node_sweep} shows that trajectory~A is contained in the admissible set and, hence, the resulting network voltages are within their limits as this dispatch trajectory is traversed. However, dispatch trajectory~B passes through the ``hole'' and results in voltage violations. Even though trajectory~A is AC admissible it requires $q_{\text{g},2}$ and $q_{\text{g},3}$ to be coordinated (i.e., stay on the trajectory) to ensure admissibility. This means that any change to one requires a change in the other and, thus, they are not considered nodal reactive capacities.
This simple example shows the need to develop methods that compute nodal reactive capacities for wind farms. This can avoid communication requirements between turbines in a wind-farm, paving the way for fast real-time control. Towards that objective, the next section develops a convex inner approximation of the non-convex \textit{DistFlow} formulation in~\eqref{eq:dist_flow}. 

\section{Convex inner approximation and optimization formulation}\label{sec:ciaOpti}

\subsection{Convex inner approximation}

In this section, we first present a compact matrix representation of the linear components~\eqref{eq:volt_rel}-\eqref{eq:reac_power_rel}.
%present the variables $P_{ij}$, $Q_{ij}$ and $v_j$ as linear functions of the power injections and the branch currents. This allows us to 
%separate the model into linear and nonlinear components and
Then, we bound the nonlinear branch current terms in~\eqref{eq:curr_rel}, $l_{ij}(P_{ij},Q_{ij},v_j)$, by a convex envelope, which leads to a convex inner approximation of~\eqref{eq:dist_flow}.

% From the \textit{incidence} matrix $B$ of the radial network and following the method adopted in~\cite{heidari2017non}, \eqref{eq:real_power_rel} and \eqref{eq:reac_power_rel} can be expressed through the following matrix equations:
% \begin{align}
%     P=p+AP-ARl \qquad Q=q+AQ-AXl, \label{eq:P_matrix} 
% \end{align}

First, define vectors $P:=[P_{ij}]_{(i,j)\in \mathcal{L}}\in \mathbb{R}^N$, $Q:=[Q_{ij}]_{(i,j)\in \mathcal{L}}\in \mathbb{R}^N$, $V:=[v_i]_{i\in \mathcal{N}}\in \mathbb{R}^N$, $p:=[p_i]_{i\in \mathcal{N}}\in \mathbb{R}^N$, $q:=[q_i]_{i \in \mathcal{N}}\in \mathbb{R}^N$, and  $l:=[l_{ij}]_{(i,j)\in \mathcal{L}}\in \mathbb{R}^N$ and matrices $R:=\text{diag}\{r_{ij}\}_{(i,j)\in \mathcal{L}}\in \mathbb{R}^{N\times N}$, $X:=\text{diag}\{x_{ij}\}_{(i,j)\in \mathcal{L}}\in \mathbb{R}^{N\times N}$, $Z^2:=\text{diag}\{z_{ij}^2\}_{(i,j)\in \mathcal{L}}\in \mathbb{R}^{N\times N}$, and $A:=[0_N \quad I_N]B-I_N$, where $I_N$ is the $N\times N$ identity matrix and $0_N$ is a column vector of $N$ rows.
Then, directly applying~\cite{heidari2017non}, we get expressions for $P$, $Q$ and $V$:
\begin{align}
    V=v_{\text{0}}\mathbf{1}_N+M_{\text{p}}p+M_{\text{q}}q-Hl,\label{eq:final_volt_rel} \\
    P=Cp-D_{\text{R}}l, \qquad Q=Cq-D_{\text{X}}l,\label{eq:P_relation}
\end{align}
where matrices $M_{\text{p}}:=2C^TRC$, \quad $M_{\text{q}}:=2C^TXC$, \quad $H:=C^T(2(RD_{\text{R}}+XD_{\text{X}})+Z^2)$ and $C:=(I_N-A)^{-1}$, $D_{\text{R}}:=(I_N-A)^{-1}AR$, and $D_{\text{X}}:=(I_N-A)^{-1}AX$ describe the network topology and impedance parameters. Note that in ~\cite{nazir2019convex}, it is proven that the matrix $(I_N-A)$ is non-singular for radial, balanced distribution networks. 
%Furthermore, the convex inner approximation in~\cite{nazir2019voltage} is valid only for purely inductive, radial, and balanced networks.
%In the current manuscript, we extend the convex formulation to any radial and balanced network, including those with mixed inductive and capacitive branches.

Clearly, \eqref{eq:final_volt_rel} and~\eqref{eq:P_relation} represent linear relationships between the nodal power injections, $(p,q)$, the branch power flows, $(P,Q)$, and node voltages $V$. 
%\st{The nonlinearity in the network is represented by~\eqref{eq:curr_rel}, as the current term $l$ is related to the power injections and node voltages in a nonlinear fashion.} 
However, setting $l=0$ and neglecting~\eqref{eq:curr_rel}, as done with the commonly used \textit{LinDist approximation}, can result in overestimating the nodal reactive capacities~\cite{nazir2019convex}. Next, we present methods for bounding the nonlinearity $l_{ij}(P_{ij},Q_{ij},v_j)$ from above and below.

%The problem being addressed in this manuscript is to determine the admissible operating range of DERs that respects the nodal voltage and line limit constraints. 
Based on the description of voltages in~\eqref{eq:final_volt_rel} and branch flows in~\eqref{eq:P_relation}, denote $l_{\text{lb}}$ and $l_{\text{ub}}$ as lower and upper bounds on $l$. Then, we can define the corresponding upper $(.)^+$ and lower $(.)^-$ bounds of $P$, $Q$ and $V$ as follows:
\begin{subequations}\label{eq:CIA_bounds}
\begin{align}
    P^+(p):=&Cp-D_{\text{R}}l_{\text{lb}} \label{eq:P_relation_1}\\
    P^-(p):=&Cp-D_{\text{R}}l_{\text{ub}} \label{eq:P_relation_2}\\
    Q^+(q):=&Cq-D_{\text{X}_+}l_{\text{lb}}-D_{\text{X}_-}l_{\text{ub}}\label{eq:Q_relation_1}\\
    Q^-(q):=&Cq-D_{\text{X}_+}l_{\text{ub}}-D_{\text{X}_-}l_{\text{lb}}\label{eq:Q_relation_2}\\
    V^+(p,q):=&v_{\text{0}}\mathbf{1}_n+M_{\text{p}}p+M_{\text{q}}q-H_+l_{\text{lb}}-H_-l_{\text{ub}}\label{eq:V_relation_1}\\
    V^-(p,q):=&v_{\text{0}}\mathbf{1}_n+M_{\text{p}}p+M_{\text{q}}q-H_+l_{\text{ub}}-H_{-} l_{\text{lb}}, \label{eq:V_relation_2}
\end{align}
\end{subequations}
 %\st{where $l_{min}$ and $l_{max}$ are the lower and upper bounds of $l$, i.e., $l_{min}\le l$ and $l_{max}\ge l$ and where}
where $D_{\text{X}_+}$ and $H_+$ include the non-negative elements of $D_{\text{X}}$ and $H$, respectively, and $D_{\text{X}_-}$ and $H_-$ are the corresponding negative elements. For example, if the network is purely inductive, then $D_{\text{X}_-}=H_-=0$ and the formulation reduces to the one presented in~\cite{nazir2019voltage}. These upper and lower bounds in~\eqref{eq:CIA_bounds} satisfy $P^-\le P \le P^+$, $Q^-\le Q\le Q^+$ and $V^-\le V\le V^+$.
%\st{However, in this work we extended the application to consider any radial distribution network through the above formulation.} 
Note that the bounds $l_\text{lb}, l_\text{ub}$ in~\eqref{eq:CIA_bounds} effectively allow us to neglect the nonlinear~\eqref{eq:curr_rel}. Thus, if we can find convex representations of these bounds, the corresponding OPF formulation will be a convex inner approximation. This is described next. 
\textcolor{black}{
\begin{remark}
While this work uses a network where the wind turbines' inductive transformers overcome the lines' capacitances, the approach presented herein readily extends to networks with arbitrary inductance/capacitance impedances as presented in~\cite{nazir2019grid}. 
\end{remark}
}
     Equation~\eqref{eq:CIA_bounds} provides a linear formulation for bounding the AC power flow equations in terms of bounds $l_\text{lb}, l_\text{ub}$ and controllable generations. %In~\cite{nazir2019convex}, we provided global, but conservative bounds on the nonlinearity based on worst case net-demand forecasts.
     %This was first presented in~\cite{nazir2019voltage}, where bounds $l_\text{lb}, l_\text{ub}$ were derived based on a nominal operating point and used to maximize voltage margins with mechanical grid assets (e.g., LTCs and capacitor-banks).
     Next, we summarize the derivation of these bounds and leverage them to formulate a novel convex inner approximation of the AC OPF to determine the nodal reactive capacities for the wind farm network.
  
Based on any nominal or predicted operating point $x^0_{ij}:=\text{col}\{P_{ij}^0, Q_{ij}^0, v_j^0\} \in\mathbb{R}^3$, the second-order Taylor series approximation for~\eqref{eq:curr_rel} can be expressed as:
\begin{align}\label{eq:T_exp}
    l_{ij} & \approx l_{ij}^0 + \mathbf{J_{ij}^\top} \mathbf{\delta_{ij}} +\frac{1}{2}\mathbf{\delta_{ij}^\top} \mathbf{H_{\text{e},ij}} \mathbf{\delta_{ij}}
\end{align}
where $l_{ij}^0 := l_{ij}(P_{ij}^0,Q_{ij}^0,v_j^0)$ are branch current flows at the operating point and $\mathbf{\delta_{ij}}(P_{ij}, Q_{ij}, v_j,x_{ij}^0)$, the Jacobian $\mathbf{J_{ij}}$ and the Hessian $\mathbf{H_{\text{e},ij}}$ are defined below:
\begin{align}
    \mathbf{\delta_{ij}}:=\begin{bmatrix} P_{ij}-P_{ij}^0\\
Q_{ij}-Q_{ij}^0\\ v_j-v_j^0\end{bmatrix} \qquad
%\begin{align}\label{eq:Jacobian}
\mathbf{J_{ij}}:=\begin{bmatrix}\frac{2P^0_{ij}}{v^0_j}\\\frac{2Q^0_{ij}}{v^0_j}\\ -\frac{(P^0_{ij})^2+(Q^0_{ij})^2}{(v^0_j)^2}\end{bmatrix}\\
    \mathbf{H_{\text{e},ij}}:=\begin{bmatrix} \frac{2}{v^0_j} && 0 && \frac{-2P^0_{ij}}{(v^0_j)^2}\\
    0 && \frac{2}{v^0_j} && \frac{-2Q^0_{ij}}{(v^0_j)^2}\\
    \frac{-2P^0_{ij}}{(v^0_j)^2} && \frac{-2Q^0_{ij}}{(v^0_j)^2} && 2\frac{(P^0_{ij})^2+(Q^0_{ij})^2}{(v^0_j)^3}
    \end{bmatrix}
\end{align}

 The expression in~\eqref{eq:T_exp} holds if we can neglect the third order term, i.e., the expression is cubic order accurate or the order of accuracy is $\mathcal{O}(||\delta||_{\infty}^3)$.
 
 Furthermore, ~\cite{nazir2019voltage} shows that $\mathbf{H_{\text{e},ij}}$ is positive semi-definite, which, together with~\eqref{eq:T_exp}, means that the lower and upper bounds of $l_{ij}$ for all $(i,j)\in \mathcal{L}$ are given by:
\begin{align}
    l_{ij}=|l_{ij}| & \approx |l_{ij}^0 + \mathbf{J_{ij}^\top}\mathbf{\delta_{ij}}+\frac{1}{2}\mathbf{\delta_{ij}^\top} \mathbf{H_{\text{e},ij}} \mathbf{\delta_{ij}}| \\
    & \le |l_{ij}^0| + |\mathbf{J_{ij}^\top}\mathbf{\delta_{ij}}|+|\frac{1}{2}\mathbf{\delta_{ij}^\top} \mathbf{H_{\text{e},ij}} \mathbf{\delta_{ij}}| \\
    & \le l_{ij}^0 + \max\{2|\mathbf{J_{ij}^\top}\mathbf{\delta_{ij}}|,|\mathbf{\delta_{ij}^\top} \mathbf{H_{\text{e},ij}} \mathbf{\delta_{ij}}|\}\label{eq:l_upper_1}\\
    \implies  l_{ij} &
     \le l_{ij}^0 + \max\{2|\mathbf{J_{ij+}}^\top\mathbf{\delta_{ij}^+}+\mathbf{J_{ij-}}^\top\mathbf{\delta_{ij}^-}|,\mathbf{\psi_{ij}}\} =: l_{\text{ub},ij} \label{eq:l_upper}\\
      l_{ij} & \ge l_{ij}^0 + \mathbf{J_{ij+}}^\top\mathbf{\delta_{ij}^-}+\mathbf{J_{ij-}}^\top\mathbf{\delta_{ij}^+} =: l_{\text{lb},ij}\label{eq:l_lower}
\end{align}
where $\mathbf{J_{ij+}}$ and $\mathbf{J_{ij-}}$ includes the positive and negative elements of $\mathbf{J_{ij}}$, $\mathbf{\delta_{ij}^+} := \mathbf{\delta_{ij}}(P_{ij}^+, Q_{ij}^+, v_j^+,x_{ij}^0)$ and $\mathbf{\delta_{ij}^-}:=\mathbf{\delta_{ij}}(P_{ij}^-, Q_{ij}^-, v_j^-,x_{ij}^0)$, and $\mathbf{\psi_{ij}}:=\max\{\mathbf{(\delta_{ij}}^{+,-})^\top \mathbf{H_{\text{e},ij}} (\mathbf{\delta_{ij}}^{+,-})\}$, which represents the largest of eight possible combinations of $P/Q/v$ terms in $\mathbf{\delta_{ij}}$ with mixed $+,-$ superscripts.
%Note that from~\eqref{eq:l_lower}, the lower bound $l_{\text{lb},ij}$ may become negative, however, we know from physics that $l_{ij}\ge 0$, which means the $l_{\text{lb},ij}$ may be overly conservative.
%{\color{blue}Also, from the expression of $\mathbf{J_{ij}}$ and $\mathbf{H_{\text{e},ij}}$, it can be seen that if $P_{ij}^0=Q_{ij}^0\approx 0$, then $l_{ij}=l_{ij}^0$ and so the lower and upper bounds do not provide information about the non-linear term $l_{ij}$.}
%To alleviate this shortcoming, Algorithm~\ref{alg1} in Section~\ref{sec:Iter_alg} presents an iterative approach that improves the nodal reactive capacity.
Thus, from~\eqref{eq:CIA_bounds},~\eqref{eq:l_upper} and~\eqref{eq:l_lower} we have a convex inner approximation of~\eqref{eq:dist_flow} that can be used to determine the nodal reactive capacities. %The convex inner approximation problem to achieve this is presented in the next section.

% The above equation can also be expressed in vector form as:
% \begin{align}\label{eq:l_upper}
%     l_{\text{max}}=\max\{l^0+2J\Delta x, \ l^0+\Delta x^T H_{\text{e}} \Delta x\}
% \end{align}

\subsection{Optimizing wind farm nodal reactive capacity}\label{sec:CIA_form} 
\begin{figure}[t]
\centering
\includegraphics[width=0.35\textwidth]{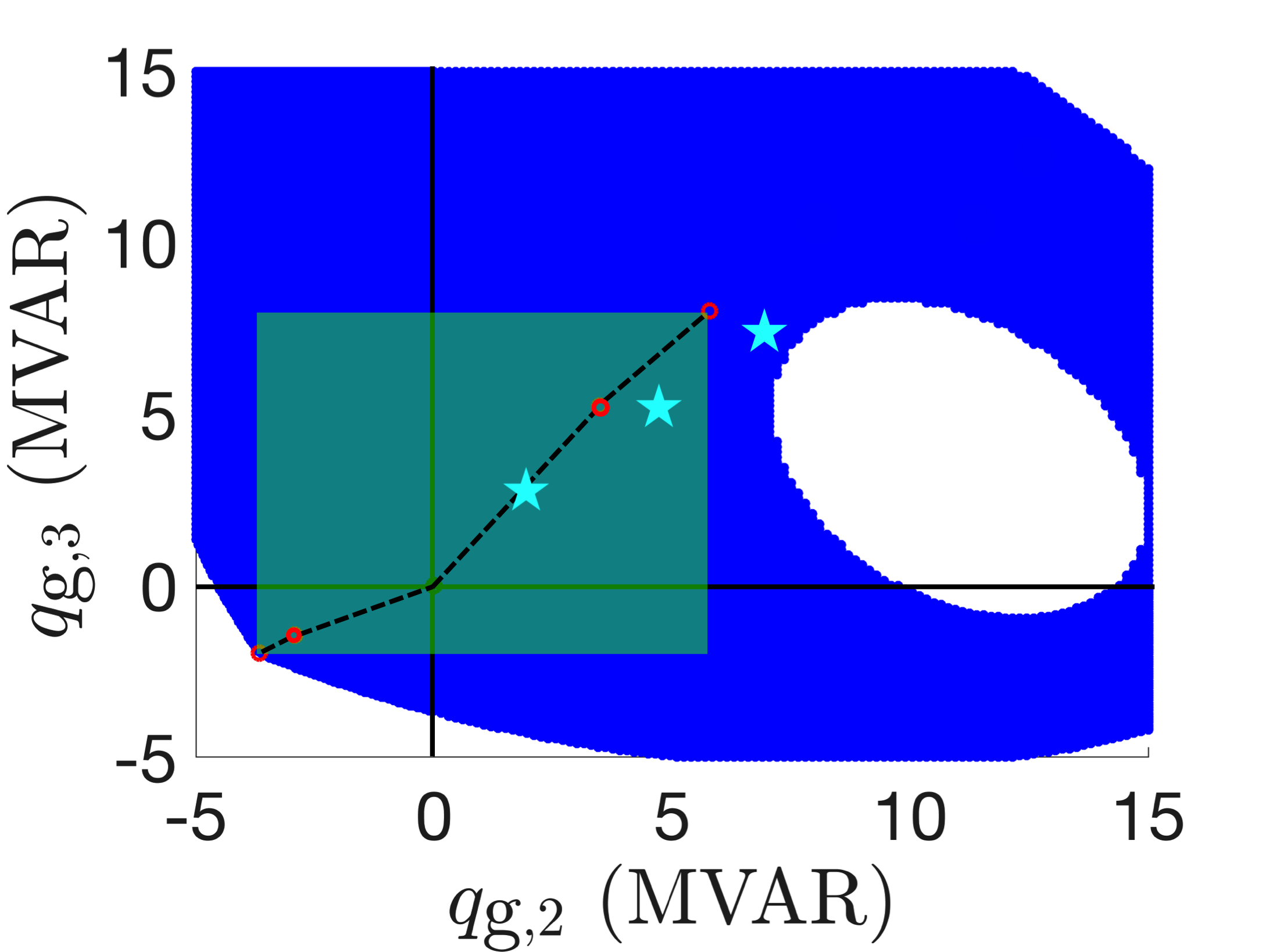}
\caption{\label{fig:3_node_inner} Algorithm is adapted to increase the admissible region via iterations (red dots), where the inner (green) set's boundary defines the nodal reactive capacities. Also, solutions obtained by the method in~\cite{Martins_2015} are indicated by the cyan stars.}
\end{figure}

The bounds from~\eqref{eq:l_upper} and~\eqref{eq:l_lower} allow us to omit~\eqref{eq:curr_rel} entirely and replace the original variables $P$, $Q$, and $V$ with their corresponding upper and lower bounds $(.)^+$ and $(.)^-$ in~\eqref{eq:CIA_bounds}. Since $(.)^+$ and $(.)^-$ are outer approximations, using them in an OPF formulation results in a feasible set that is contained in the original, non-convex AC OPF, which means that (P1) and (P2) below represent convex inner approximations and can be used to determine the wind farm nodal reactive capacity:
%Based on these bounds, we can now formulate the complete convex inner approximation OPF problem that determines the nodal reactive capacity across the feeder shown in (P1) as:
%\begin{subequations}\label{eq:P1}
\begin{align}
\text{(P1)}  \quad q^{+}=\arg \min_{q_i} \  - Q_\text{10}^{-}(q)+\sum_{i=1}^{N}f_i(q_{i})&\\
\text{s.t.}  \quad  \eqref{eq:P_relation_1}-\eqref{eq:V_relation_2}, \eqref{eq:l_upper}, \eqref{eq:l_lower}&\\
    \underline{V}\le     V^-(q) \quad V^+(q)\le &\overline{V}   \label{eq:P1_V_a}\\
    % (P^+(p,q))^2+  (Q^+(p,q))^2 \le & \overline{l}\underline{V} \label{eq:P1_L1_limit}\\
    % (P^+(p,q))^2+  (Q^-(p,q))^2 \le & \overline{l}\underline{V}  \\
    % (P^-(p,q))^2+ (Q^+(p,q))^2 \le & \overline{l}\underline{V} \\ (P^-(p,q))^2+ (Q^-(p,q))^2 \le & \overline{l}\underline{V}\label{eq:P1_L2_limit}
    l_{\text{ub}}\le  \overline{l} \quad \underline{q}\le q\le \overline{q} \label{eq:P1_lmax_a}
\end{align}
\begin{align}
\text{(P2)}  \quad q^{-}=\arg \min_{q_i} \  Q_\text{10}^{+}(q)+\sum_{i=1}^{N}f_i(q_{i})&\\
\text{s.t.}  \quad  \eqref{eq:P_relation_1}-\eqref{eq:V_relation_2}, \eqref{eq:l_upper}, \eqref{eq:l_lower}&\\
    \underline{V}\le     V^-(q) \quad V^+(q)\le &\overline{V}   \label{eq:P1_V_b}\\
    % (P^+(p,q))^2+  (Q^+(p,q))^2 \le & \overline{l}\underline{V} \label{eq:P1_L1_limit}\\
    % (P^+(p,q))^2+  (Q^-(p,q))^2 \le & \overline{l}\underline{V}  \\
    % (P^-(p,q))^2+ (Q^+(p,q))^2 \le & \overline{l}\underline{V} \\ (P^-(p,q))^2+ (Q^-(p,q))^2 \le & \overline{l}\underline{V}\label{eq:P1_L2_limit}
    l_{\text{ub}}\le  \overline{l} \quad \underline{q}\le q\le \overline{q} \label{eq:P1_lmax_b}
\end{align}
%\end{subequations}
%\st{where $M \in \mathbb{R}_+$ is a large number and $z_i\in \mathbb{Z}$ is binary. The "big M" method is used here to represent \eqref{eq:l_lower} in a convex form.} 
where $Q_{10}^-$ and $Q_{10}^+$ are the lower and upper bound of the reactive power flow in line connecting nodes $1$ and $0$ (also called $Q_{\text{head}}$),~\eqref{eq:P1_V_a} and~\eqref{eq:P1_lmax_a} ,\eqref{eq:P1_V_b} and~\eqref{eq:P1_lmax_b}, ensure that any feasible dispatch $q$ from (P1) and (P2) satisfies nodal voltages and branch flows in the original AC OPF based on~\eqref{eq:dist_flow}. To determine the nodal reactive capacity, we must solve (P1) for the upper capacity $q^+$ and (P2) for the lower capacity $q^-$. Thus, the objective function components, $f_i(q_i)$, must be designed to engender $q_{i}^-$ and $q_{i}^+$. For example, when computing $q_{i}^-$, we can choose $f_i(q_{i}):=\alpha_i q_{i}\,$ and, for $q_{i}^+$, we can designate $f_i(q_{i}):=-\alpha_i q_{i}$, where $\alpha_i$ is the relative priority  of nodal reactive capacity at node $i$. Clearly, the objective function determines how nodal reactive capacities are allocated over the network, but objective function design will be explored in future work and is outside the scope of this paper. %, e.g., choosing objective function such as $\pm \alpha_i\log(q_{i})$ could result in a different allocation of nodal reactive capacity over the nodes as compared with $\pm \alpha_i q_{i}$ . %The design of the objective function represents an interesting future extension into incentive mechanism and rate design~\cite{perez2017regulatory}.
\textcolor{black}{
\begin{remark}
In this work, we use simple box constraints for the nodal reactive power limits. However, other convex constraints, such as quadratic reactive capability constraints (e..g, $p_i^2 + q_i^2 \le (S_i^\text{max})^2$ for apparent power limit, $S_i^\text{max}$) can readily be included in the formulation. For example, an extensive comparison of different active and reactive power inverter schemes (and resulting constraints) is presented within the CIA formulation in~\cite{nazir2019grid}.
\end{remark}
}
The optimization problem in (P1) and (P2) is applied to the 3-node example shown in Fig.~\ref{fig:3_node_model} to determine the inner convex set. The results for this example network are shown in Fig.~\ref{fig:3_node_inner}, where the green rectangular set is the inner approximation obtained through (P1) and (P2). In this example, we also adapt the approach in~\cite{nazir2019grid} to reactive power to iteratively expand the nodal reactive capacities as indicated by the red dots in Fig.~\ref{fig:3_node_inner}. A detailed description of the iterative approach can be found in~\cite{nazir2019grid} and is beyond the scope of this paper. The convex set in Fig.~\ref{fig:3_node_inner} allows for dispatching reactive power resources without need for coordination among them, while at the same time guaranteeing satisfaction of network constraints. For comparison we also apply the method from~\cite{Martins_2015} to this example network. Points resulting from this method are shown by the cyan colored stars in Fig.~\ref{fig:3_node_inner}.
%While (P1) ensures AC admissibility at the nodal reactive capacity values, it is natural to consider what happens when the nodal flexibility is below the rated capacity. That is, are all injections within the nodal reactive capacity range guaranteed to be admissible across all the nodes? \cite{nazir2019grid} answers this question by providing analytical guarantees of admissibility for the nodal hosting capacity, $\Delta q_{\text{g}}$, and then presents an iterative algorithm to successively improve $\Delta q_{\text{g}}$.

\subsection{Simulation resulting on nodal reactive capacity}

\begin{figure}[t]
\centering
\includegraphics[width=0.4\textwidth]{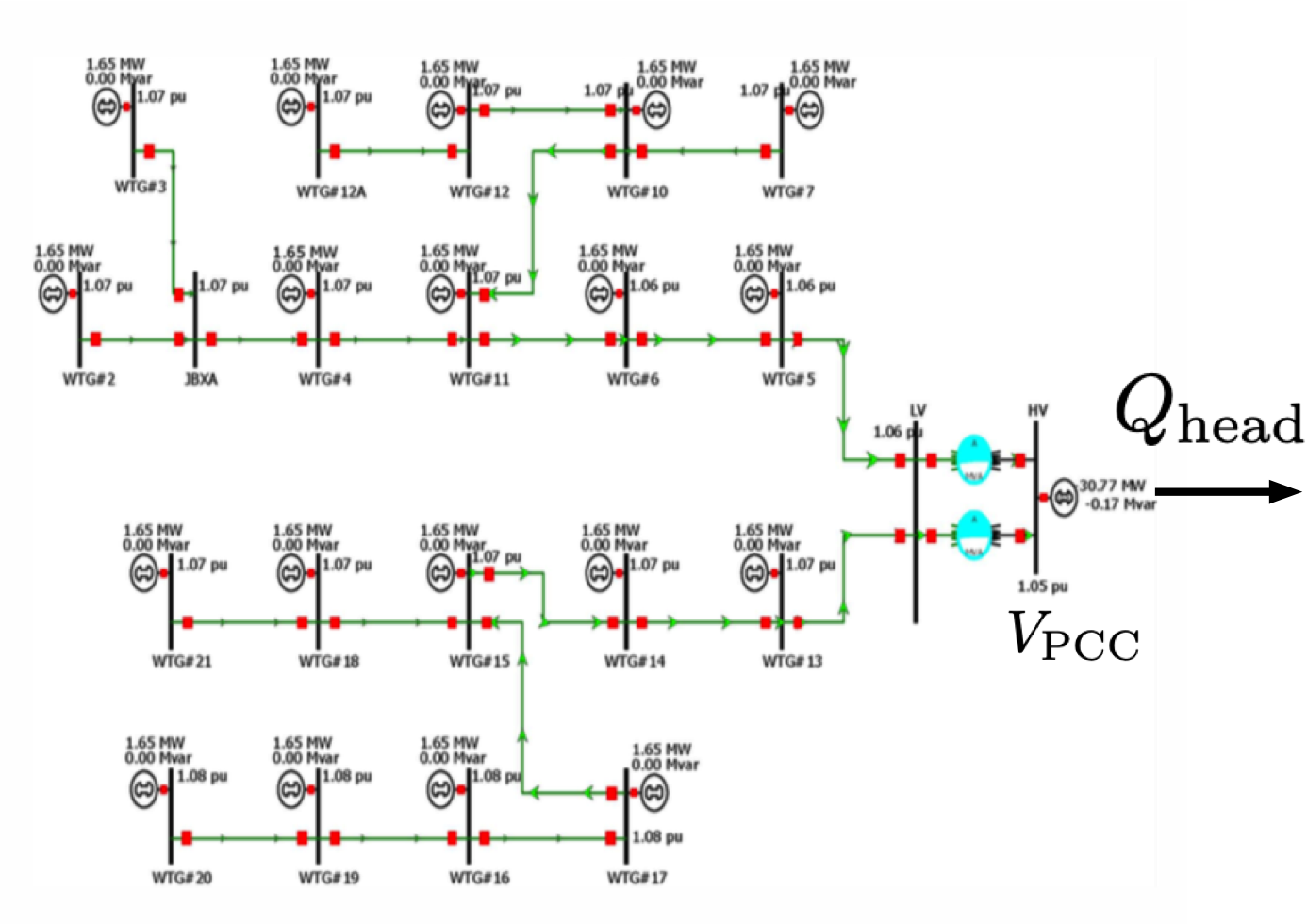}
\caption{\label{fig:wind_farm_network}Layout of the 19 turbine wind farm.}
\end{figure}

To showcase the effectiveness of the proposed convex inner approximation method, we conduct simulation-based analysis on the 19-turbine wind farm in Fig.~\ref{fig:wind_farm_network}. We determine the reactive power capability of the wind farm and the nodal reactive capacities for each turbine. We consider three scenarios similar to those in~\cite{Martins_2015}, determine the resulting wind farm reactive power capabilities and compare the results with~\cite{Martins_2015}. These three scenarios compare the reactive power capability of the wind-farm under different active power conditions. For reference, the method in~\cite{Martins_2015} broadcasts a reactive power reference signal to each wind turbine. The local control scheme of each turbine seeks to maintain the reactive power output at that broadcast reference value. However, turbine protection will override reactive power control, if necessary, to ensure the turbine's terminal voltage does not deviate beyond its limits.

We will also compare the results obtained through convex inner approximation with those obtained using a convex relaxation method (outer approximation) and by solving the full non-linear model (local solution). Each of the 19 turbines in the network is rated at 1.65 MW and can operate in a voltage range of 0.9 to 1.1 pu. Each of the turbines has reactive power capability of $[-0.5,0.5]$~MVAr. The scenario descriptions are listed below:
\begin{itemize}
    \item \textbf{Scenario~1}: All wind turbines are operating at their active power generation capacity of $p_j=1.65$ MW.
    \item \textbf{Scenario~2}: The active power generation is half the nameplate capacity of the wind farm, with the turbines at the ends of branches operating at full capacity $p_j=1.65$~MW, and the others producing zero active power.
    \item \textbf{Scenario~3}: The active power generation is half the nameplate capacity, with the turbines operating at full and zero capacity swapped from scenario~2.
\end{itemize}

\begin{figure*}[h]
\subfloat
[\label{fig:scenario_1}]{
\includegraphics[width=0.3\linewidth]{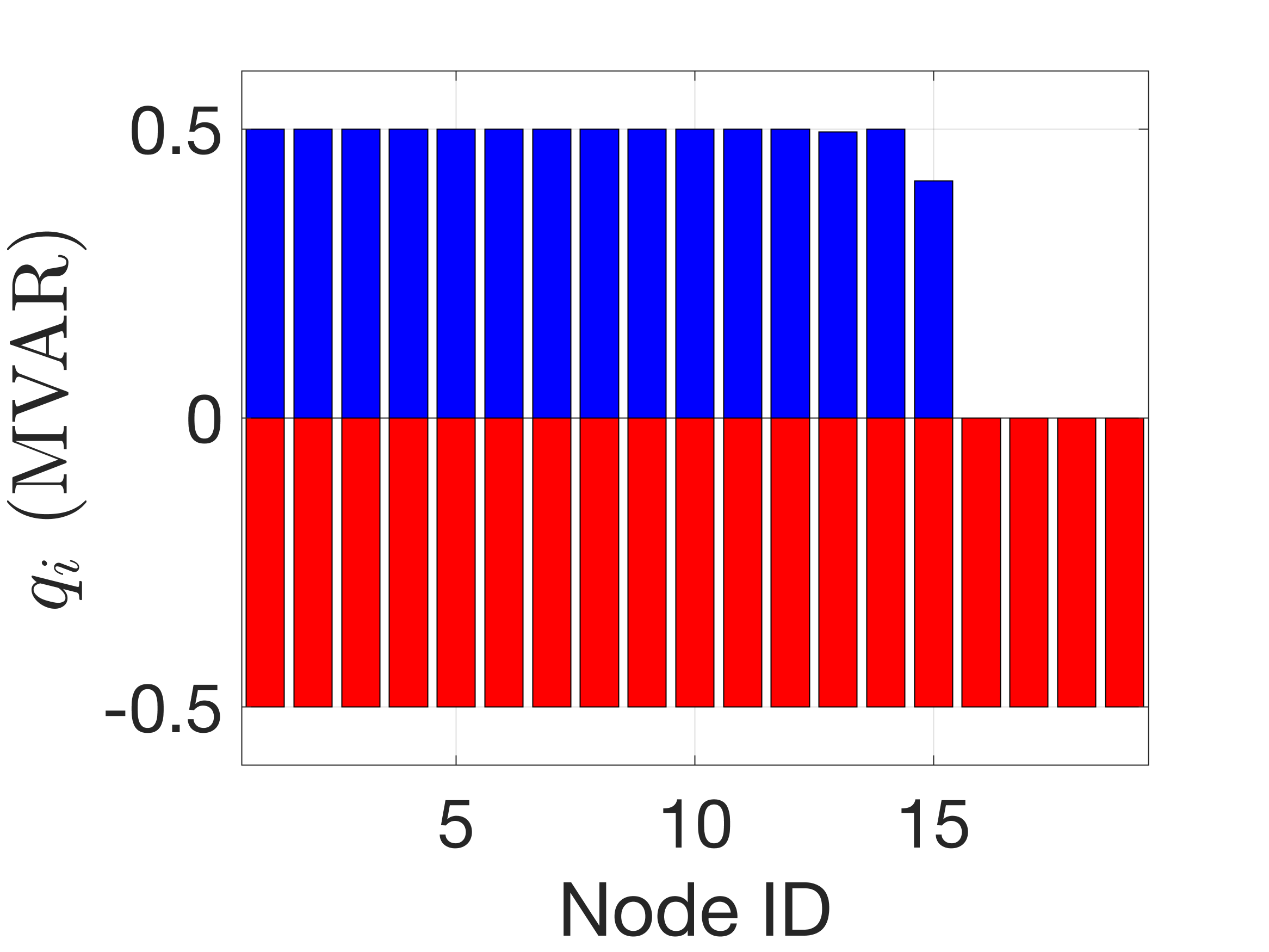}}
\hfill
\subfloat
[\label{fig:scenario_2}]{\includegraphics[width=0.3\linewidth]{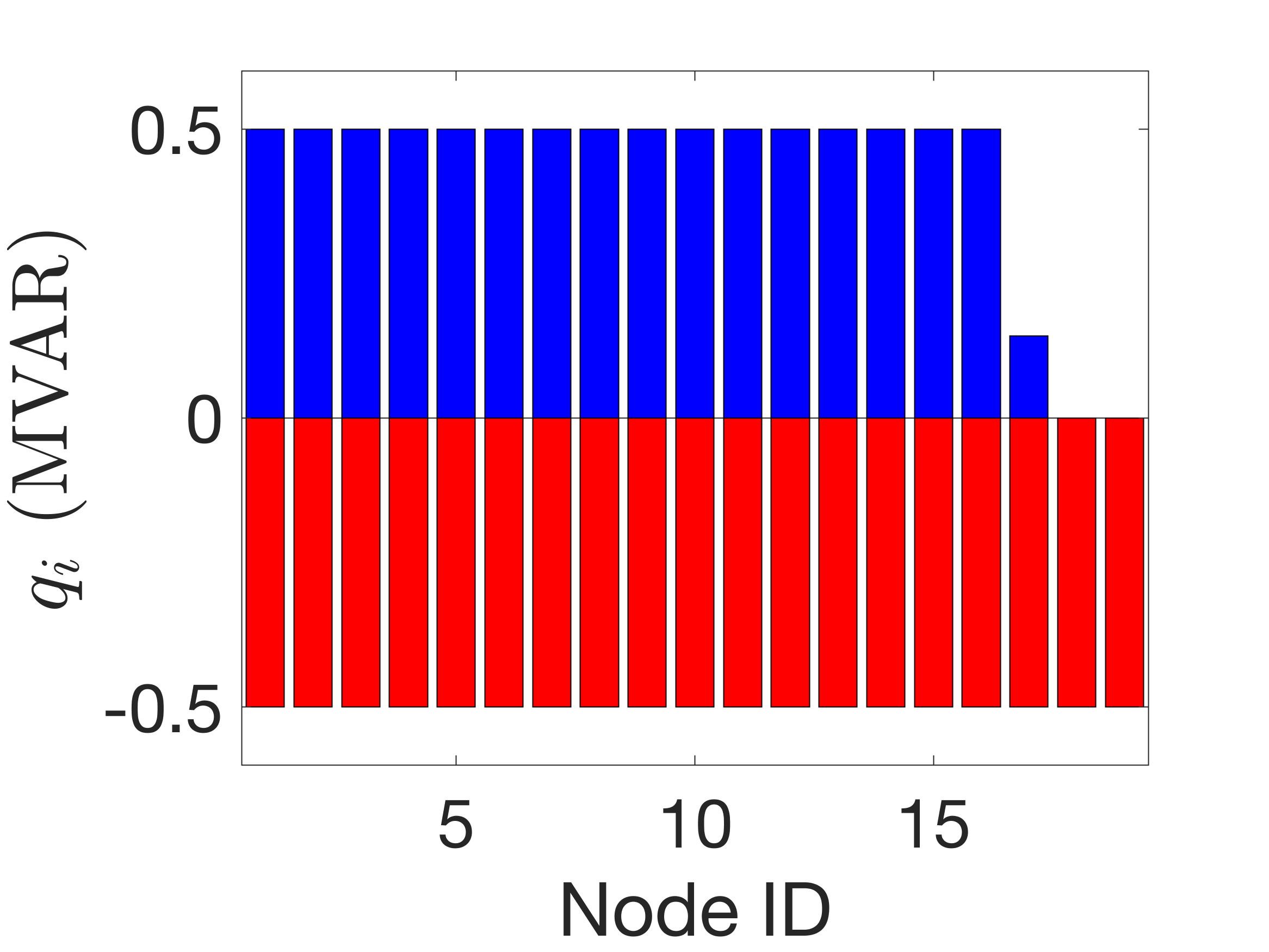}}
\hfill
\subfloat[\label{fig:scenario_3}]{\includegraphics[width=0.3\linewidth]{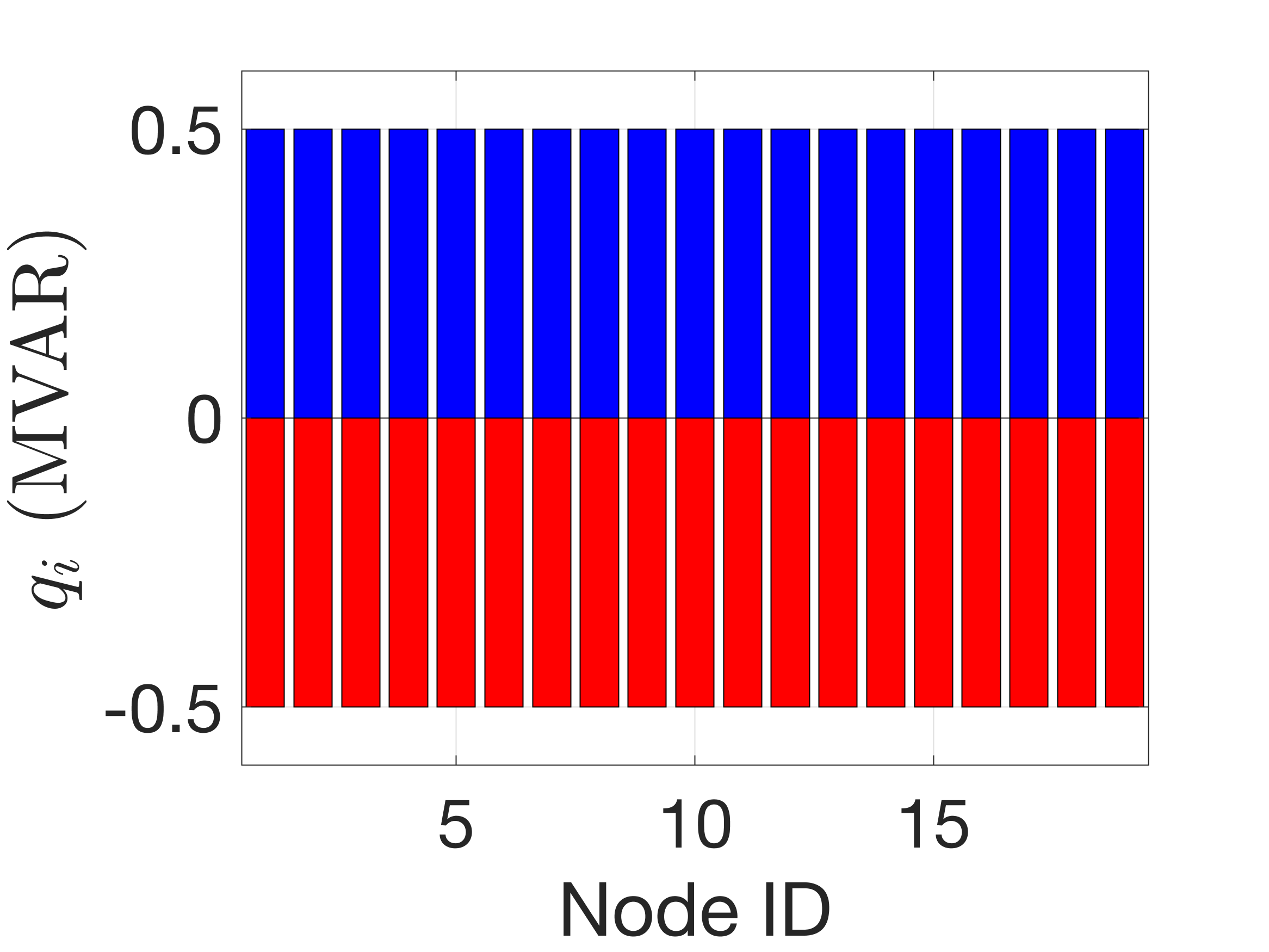}}
\caption{Nodal reactive capacities for (a) scenario~1, (b) scenario~2 and (c) scenario~3.}% In (a), all turbines $i$ operate at full generating capacity (i.e., $p_i=1.65$ MW).  For (b), scenario~2 considers the wind farm operating at half of its nameplate capacity with turbines at leaf nodes  operating at full capacity (i.e., $p_i = 1.65$ MW) and those closer to the PCC not producing any active power (i.e., $p_i = 0$ MW). In (c), the wind farm also operates at half of nameplate capacity, however, in this case, the active power generation is flipped from Scenario~2 between turbines at leaf nodes and near the head-node.}
\end{figure*}

Applying (P1) and (P2) to scenarios 1, 2, and 3 gives the nodal reactive capacities shown in Figs.~\ref{fig:scenario_1},~\ref{fig:scenario_2} and~\ref{fig:scenario_3}, respectively. Specifically, Figs.~\ref{fig:scenario_1} and~\ref{fig:scenario_2} show that the reactive power capacity is limited in scenarios~1 and~2. This is due to the upper nodal voltage limits at the turbines, which constrain the amount of reactive power support that can be provided by the entire wind farm. In contrast, no voltage limits are binding in scenario~3, allowing full reactive power support. %Scenarios~2 and~3 provide more reactive power support since the wind turbines produce less active power, allowing greater reactive power generation. 
To provide insights into the conservativeness of the CIA method, we compare the total reactive power capability of the wind farm network with results from three other approaches and formulations, the method presented in~\cite{Martins_2015}, a convex relaxation (from~\cite{gan2014exact}), and the NLP based on~\eqref{eq:dist_flow}. {\color{black} The comparison is provided in Table~\ref{table_comp} and indicates that the proposed CIA method can compute practically relevant capacities (i.e., CIA does not result in overly conservative bounds) across three different active power (wind generation) scenarios and against three different approaches/formulations. Specifically, despite being an inner approximation (i.e., inherently conservative), the CIA method's iterative approach results in reactive capacities (i.e., $Q_\text{head}$) that compare well overall against both a convex relaxation method, which is an outer approximation (i.e., overly optimistic), and the non-convex NLP-based method.}

\begin{figure}[h]
\centering
\includegraphics[width=0.48\textwidth]{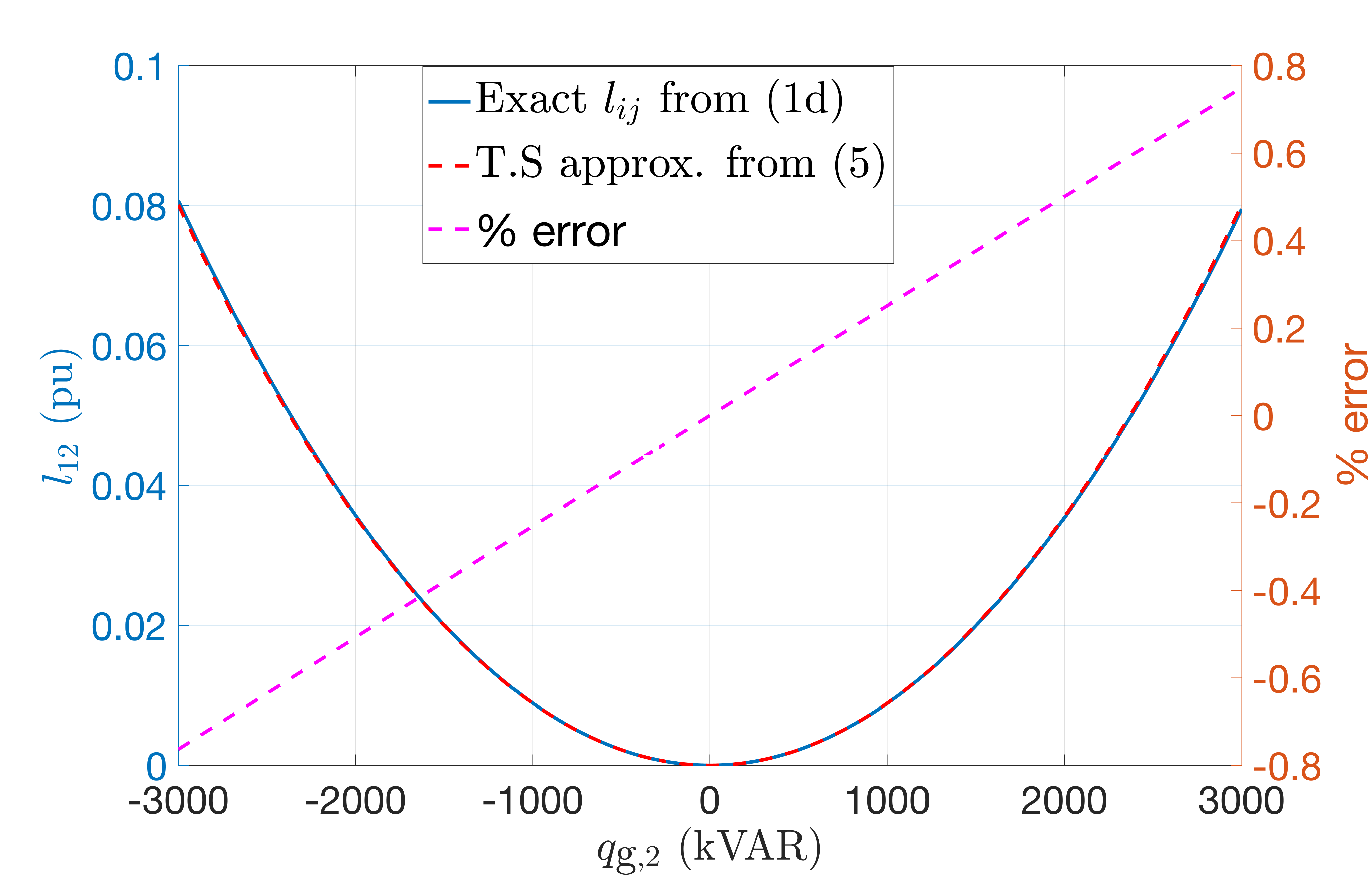}
\caption{\label{fig:TS_approx_reac} {\color{black} Comparison of the second-order approximation from~\eqref{eq:T_exp} with the original (nonlinear) expression of $l$ in~\eqref{eq:curr_rel} for the 19-turbine wind farm. The maximum absolute approximation error is found at branch $l_{12}$ (head branch) when the reactive net-injection at node~$2$ ($q_{\text{g},2}$) is $-3000$kVAR and is $6.1\times 10^{-4}$pu and the relative error is always within $0.8\%$.}}
\end{figure}

\textcolor{black}{To further highlight that the CIA method is not overly conservative, Fig.~\ref{fig:TS_approx_reac} compares the accuracy of the second-order approximation in~\eqref{eq:T_exp} with the non-linear expression from~\eqref{eq:curr_rel} for the 19-turbine wind farm in Fig.~\ref{fig:wind_farm_network}. Note that, the worst-case approximation error for $l_{ij}$ is less than $6.1\times 10^{-4}$pu over the wide range of net injections $[-3000,3000]$ kVAR.} 
 %Thus, we can assume that the expression in~\eqref{eq:T_exp} is sufficiently accurate and omit the higher-order terms.

\begin{remark}
%In comparing the performance of the CIA method with the method from~\cite{Martins_2015}, please consider the following. 
The method in~\cite{Martins_2015} achieves a larger reactive power capacity $Q_\text{head}$ than the CIA method in part because it does not restrict the reactive injections to be positive at all nodes. (This requirement is enforced in the CIA approach to obtain nodal capacities, $q_i^+$.) This additional degree of freedom enlarges the operating range. From a practical perspective, the two methods have different implementations. While the CIA method provides \textit{a~priori} predictions of $Q_\text{head}$ and corresponding decoupled, network-aware operating ranges, $[q_i^-, q_i^+]$, it requires full knowledge of the wind farm (network parameters and up-to-date active power generation) and centralized computing to determine and broadcast local reactive power nodal capacities. On the other hand, the method in~\cite{Martins_2015} requires only local sensing and control at each wind turbine, but cannot predict the wind farm's total $Q_\text{head}$ capability in advance.
\end{remark}

\begin{table}[!t]
\caption{Comparison of $Q_\text{head}$ for different schemes}
\label{table_comp}
\centering
\begin{tabular}{rlll}
\toprule
{Scheme} & {Scenario~1} & {Scenario~2} & {Scenario~3} \\
\midrule
CIA-based & [-9.9,7.0] & [-9.8,7.9] & [-9.6,9.4]\\
Nonlinear & [-9.9,7] & [-9.8,7.9] & [-9.6,9.4]\\
Relaxation & [-9.9,9.1] & [-9.8,9.3] & [-9.6,9.4]\\
From~\cite{Martins_2015} & $\le 7.3$ & $\le 8.0$ & $\le 9.4$\\             
\bottomrule
\end{tabular}
\end{table}

\begin{figure}[t]
\centering
\includegraphics[width=0.3\textwidth]{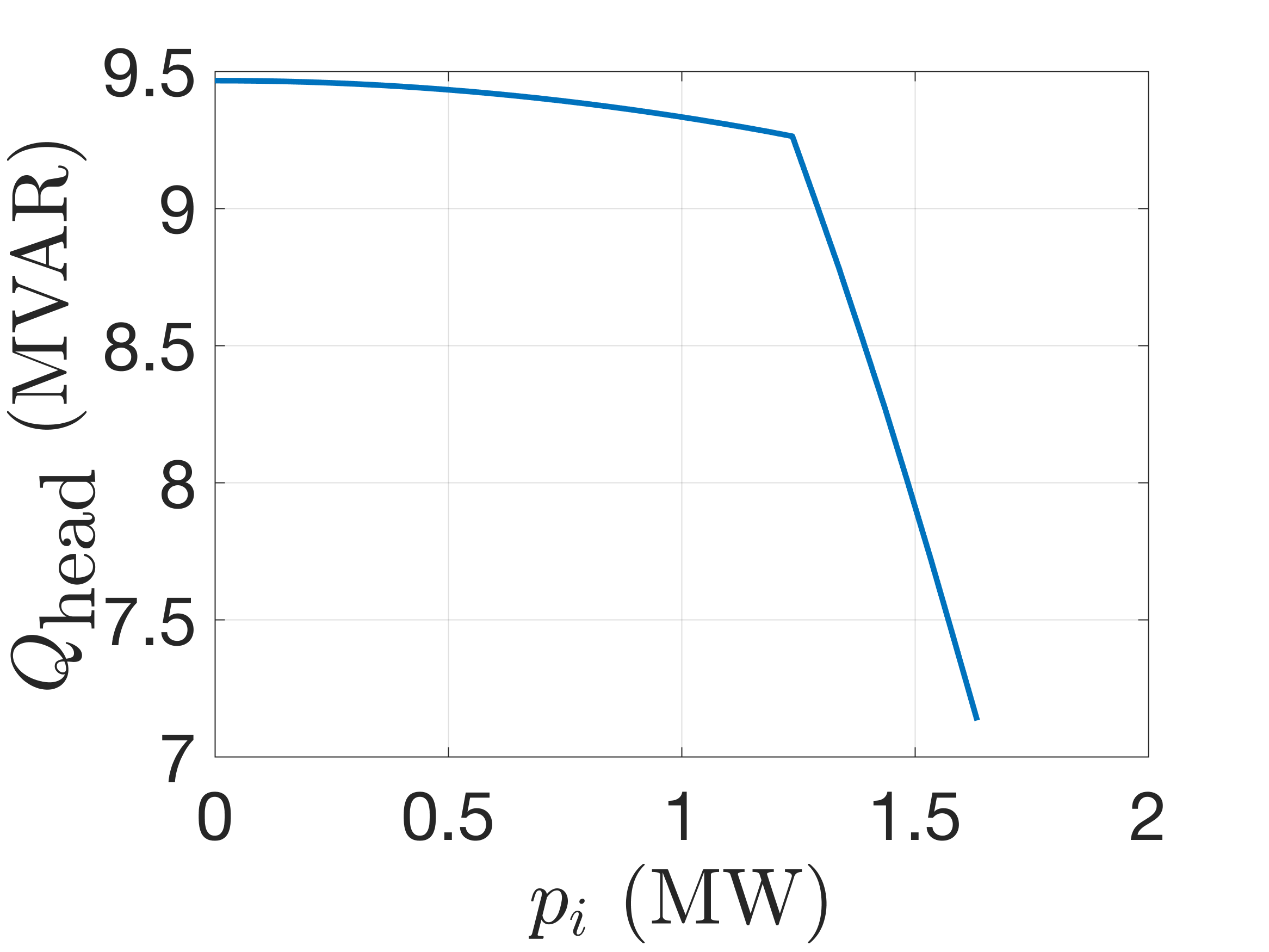}
\caption{\label{fig:P_vs_Q}Change in wind farm reactive power capacity with change in active power generation.}
\end{figure}

These scenarios indicate that the reactive power capacity of the wind farm changes as active power generation changes. Specifically, there are two constraints that are responsible for limiting the reactive power capacity: 1)  the (local) voltage constraints, and 2) the (local) turbine reactive power limit. To understand the effects of these constraints on the aggregate wind farm reactive power capacity, Fig.~\ref{fig:P_vs_Q} shows the relationship between active power generation $p_i$ (with all turbines producing the same $p_i \in [0,1.65]$ MW) and the corresponding reactive power capacity of the entire wind farm at the head-node, $Q_{\text{head}}$. As the active power generation increases from $0$ (i.e., no turbines producing active power), the wind farm's reactive power capacity decreases slowly due to increasing losses. In this initial phase, the reactive power is limited only by the upper reactive power limits $\overline{q}_i$ of the turbines. However, at around $p_i=1.2$~MW, some of the network and turbine voltage constraints become active causing the reactive power capacity to reduce sharply. %In this simulation, the active power generation is increased uniformly across all generators from 0 to rated capacity of 1.65 MW.

\section{Real-time voltage regulation algorithm}\label{sec:rtAVRdisagg}

Based on the computed nodal reactive capacities, a real-time control algorithm can be developed to dispatch reactive-power set-points to each turbine in the wind farm in order to provide voltage support at the PCC\@. We again consider the 19-turbine wind farm network shown in Fig.~\ref{fig:wind_farm_network}, where all the wind turbines are operating at full capacity of 1.65~MW\@. It is assumed that the wind farm is connected to a large power system at the PCC and the high-voltage transmission system is represented by a Th\'evenin equivalent circuit with a voltage source denoted by the voltage $V_{\text{grid}}$. The aim of this controller is to maintain the wind farm PCC voltage $V_{\text{PCC}}$ at $1.06$~pu in the presence of power system disturbances. It is assumed that the grid experiences a fault, which causes $V_{\text{grid}}$ to suddenly drop. The grid's voltage then recovers slowly as part of a recovery event, which is modeled as a ramp.

\begin{figure}[t]
\includegraphics[width=0.5\textwidth]{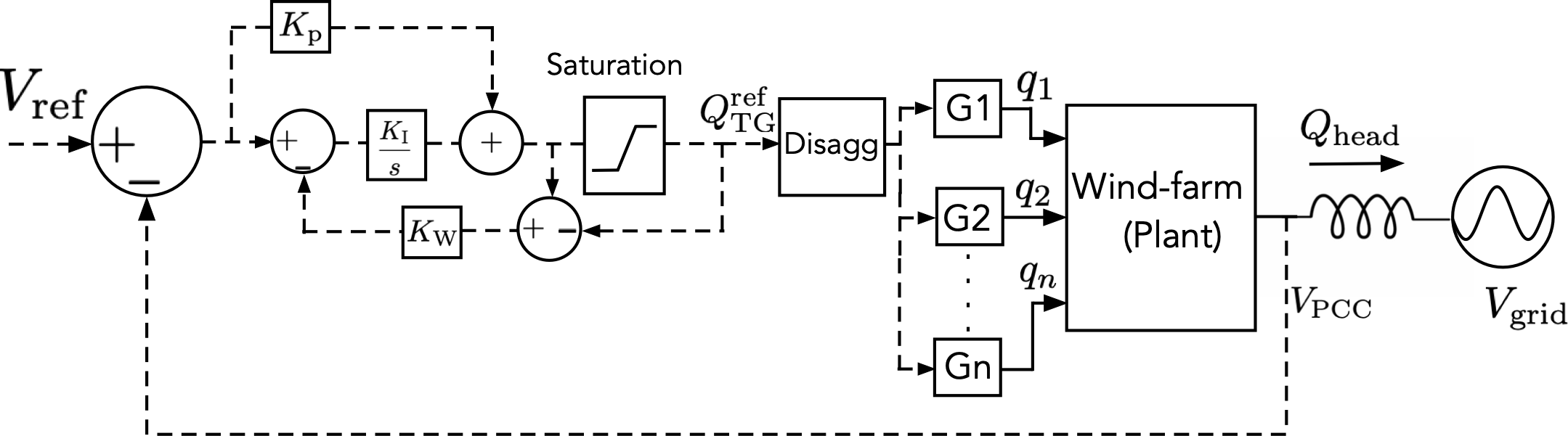}
\caption{\label{fig:windfarm_control} Proposed control scheme for real-time disaggregation and grid voltage support.}
\end{figure}

The feedback control scheme is shown in Fig.~\ref{fig:windfarm_control}. A PI controller is used to regulate the reactive power dispatch based on the deviation of $V_{\text{PCC}}$ from the reference value $V_{\text{ref}}$. The saturation block restricts the total wind farm reactive power reference $Q_{\text{TG}}^{\text{ref}}$ to within the pre-calculated maximum and minimum wind turbine reactive power limits. A standard anti-windup mechanism is implemented to ensures that if $Q_{\text{TG}}^{\text{ref}}$ saturates then the PI controller's integrator does not wind up. The `Disagg' block disaggregates the reactive power reference $Q_{\text{TG}}^{\text{ref}}[k]$ at time-step $k$ amongst the individual turbines relative to their nodal reactive capacities (i.e., $q_i^-$ and $q_i^+$) as
\begin{align}
    q_i[k]=
    \begin{cases}
    \frac{q_i^+}{\sum_iq_i^+}Q_{\text{TG}}^{\text{ref}}[k] \qquad Q_{\text{TG}}^{\text{ref}}[k]>0\\
    \frac{q_i^-}{\sum_iq_i^-}Q_{\text{TG}}^{\text{ref}}[k] \qquad Q_{\text{TG}}^{\text{ref}}[k]<0.
    \end{cases}
\end{align}
The PCC voltage $V_{\text{PCC}}$ is measured and fed back to achieve closed-loop tracking of the voltage reference $V_{\text{ref}}$.

We compare this disaggregation scheme with a grid-agnostic scheme and contrast the results of the real-time controller under these two approaches. The disaggregation schemes can be summarized as:
\begin{itemize}
    \item \textbf{CIA-based disaggregation}: Disaggregation is proportional to the nodal reactive capacities ($q^-,q^+$) and the saturation block uses the computed reactive power capacity limits for the wind farm.
    \item \textbf{Grid-agnostic disaggregation}: Disaggregation is proportional to the turbine reactive power capacities ($\underline{q},\overline{q}$) and the saturation block uses the sum of those turbine capacity values. This scheme ignores the collector network constraints.
    %\item \textbf{Disaggregation based on~\cite{Martins_2015}}: In this scheme we employ the nodal capacities determined in~\cite{Martins_2015} and disaggregate proportionally based on them. The saturation block is also determined based on the nodal capacities calculated from~\cite{Martins_2015}.
\end{itemize}

\begin{figure*}[t]
\subfloat
[\label{fig:Vref_grid_agnostic}]{
\includegraphics[width=0.32\linewidth]{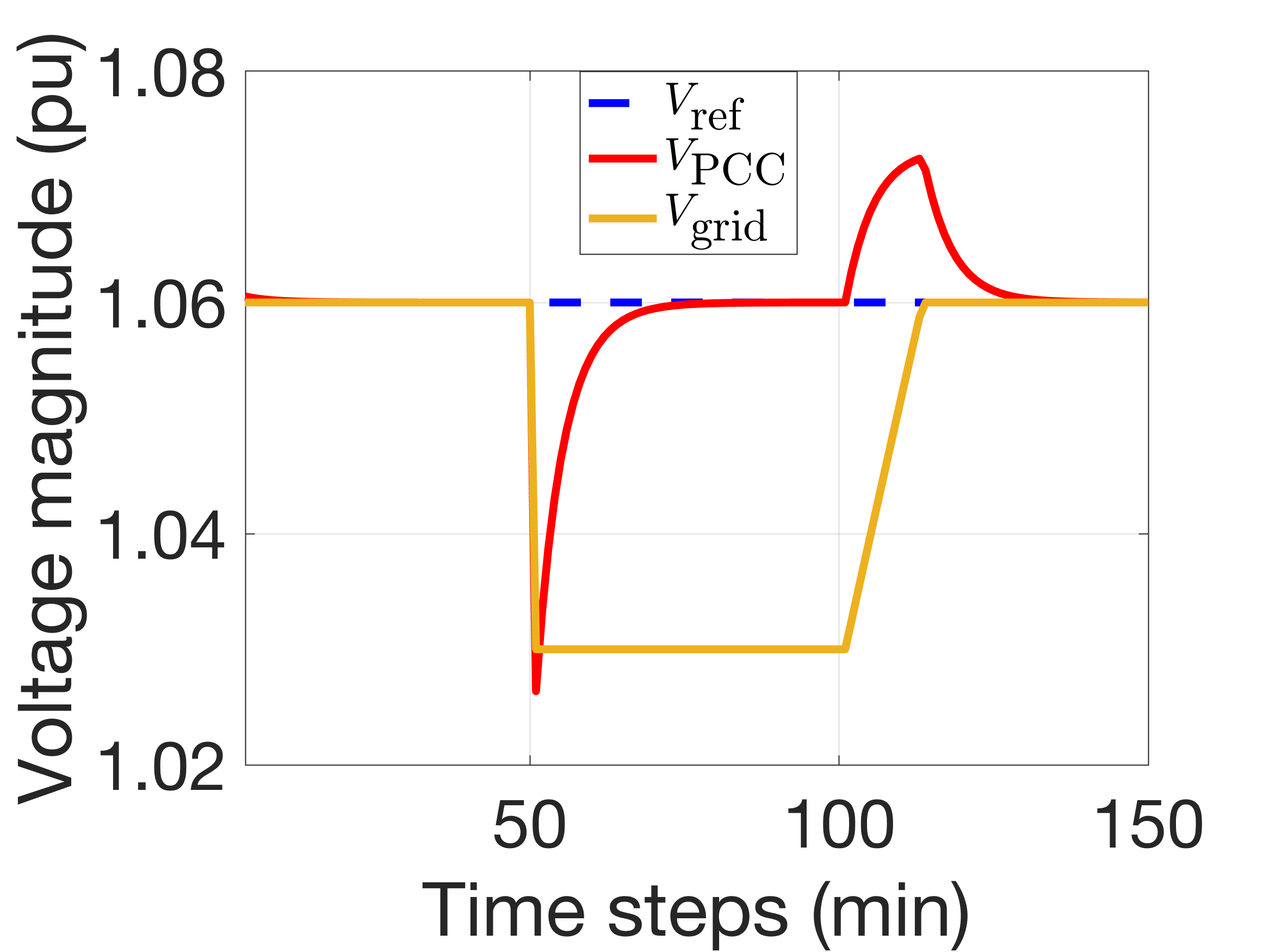}}
\hfill
\subfloat
[\label{fig:Qref_grid_agnostic}]{
\includegraphics[width=0.32\linewidth]{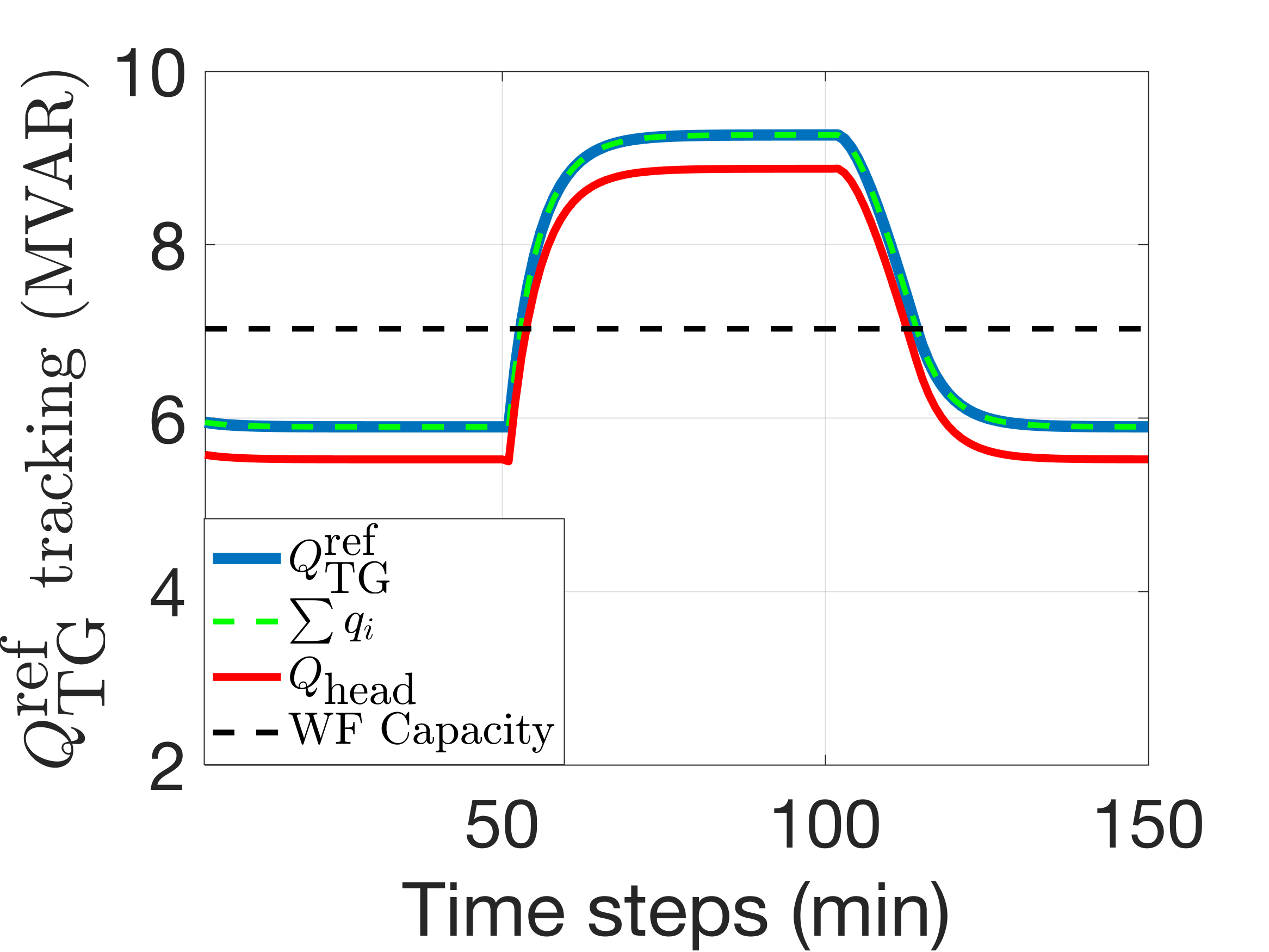}}
\hfill
\subfloat
[\label{fig:Vwf_grid_agnostic}]{
\includegraphics[width=0.32\linewidth]{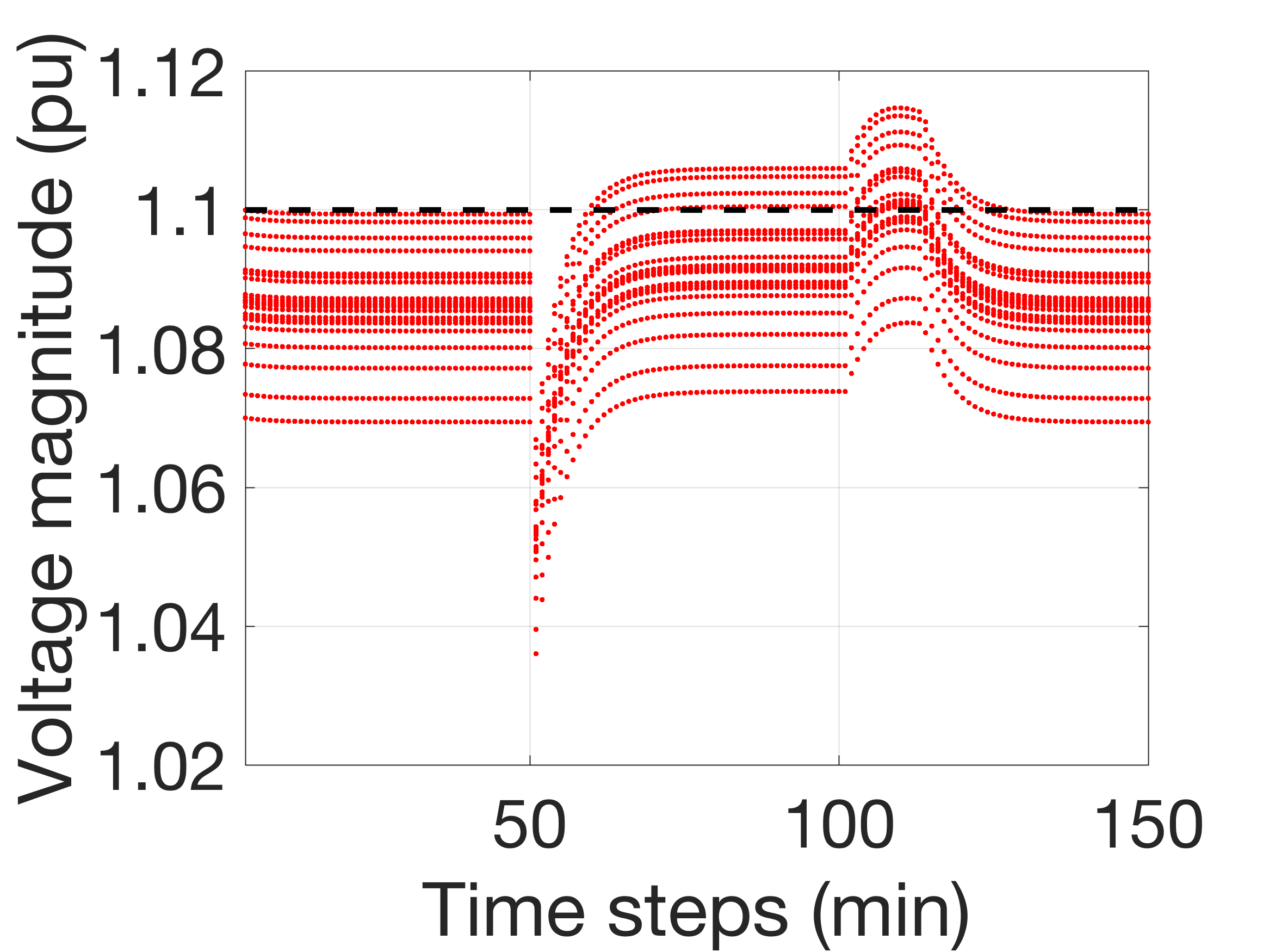}}
\\
% \subfloat
% [\label{fig:Vref_Martins}]{\includegraphics[width=0.33\linewidth]{images/V_ref_head_Martins.png}}
% \hfill
\subfloat[\label{fig:Vref_CIA}]{\includegraphics[width=0.33\linewidth]{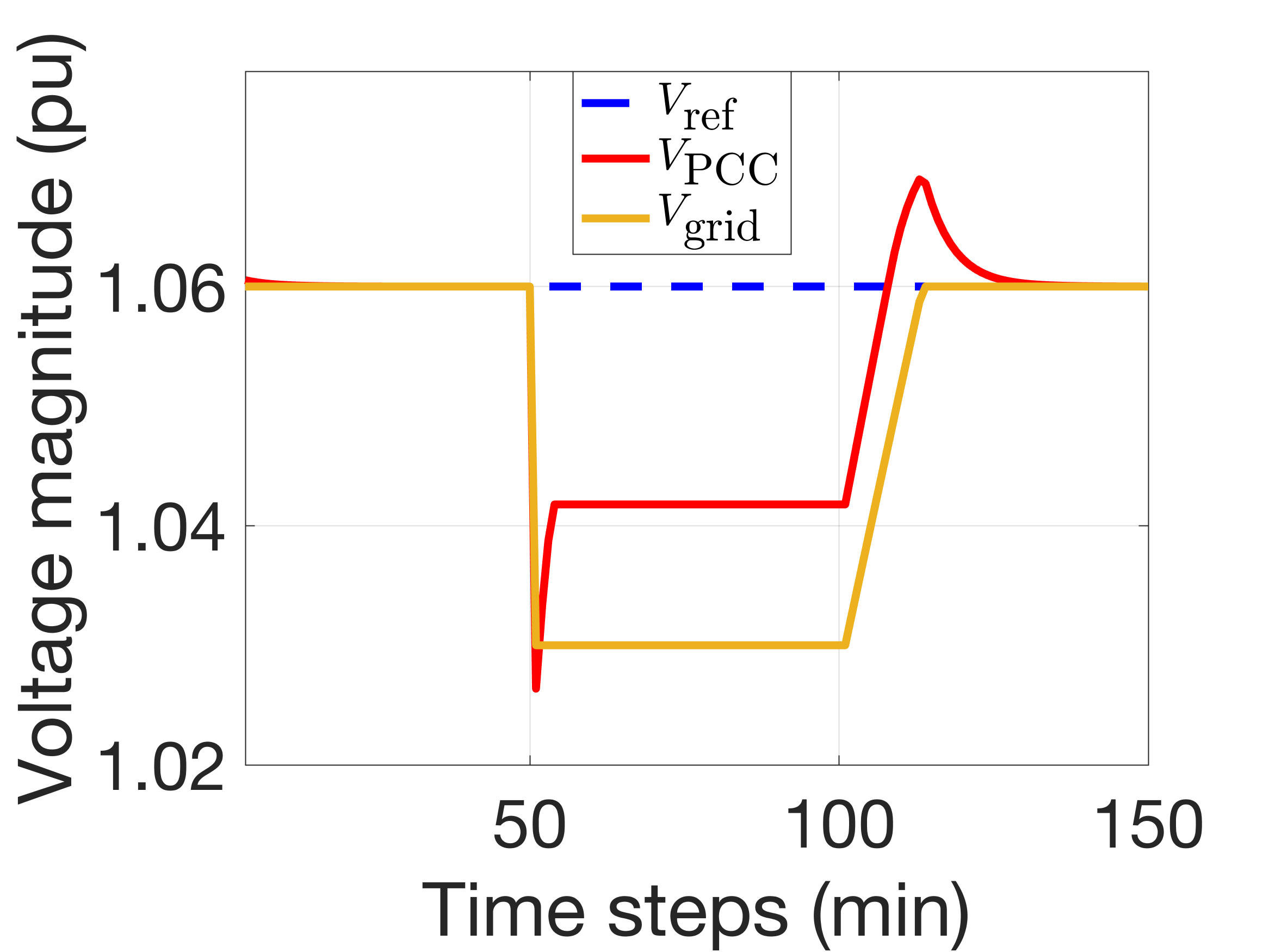}}
\hfill
\subfloat[\label{fig:Qref_CIA}]{\includegraphics[width=0.33\linewidth]{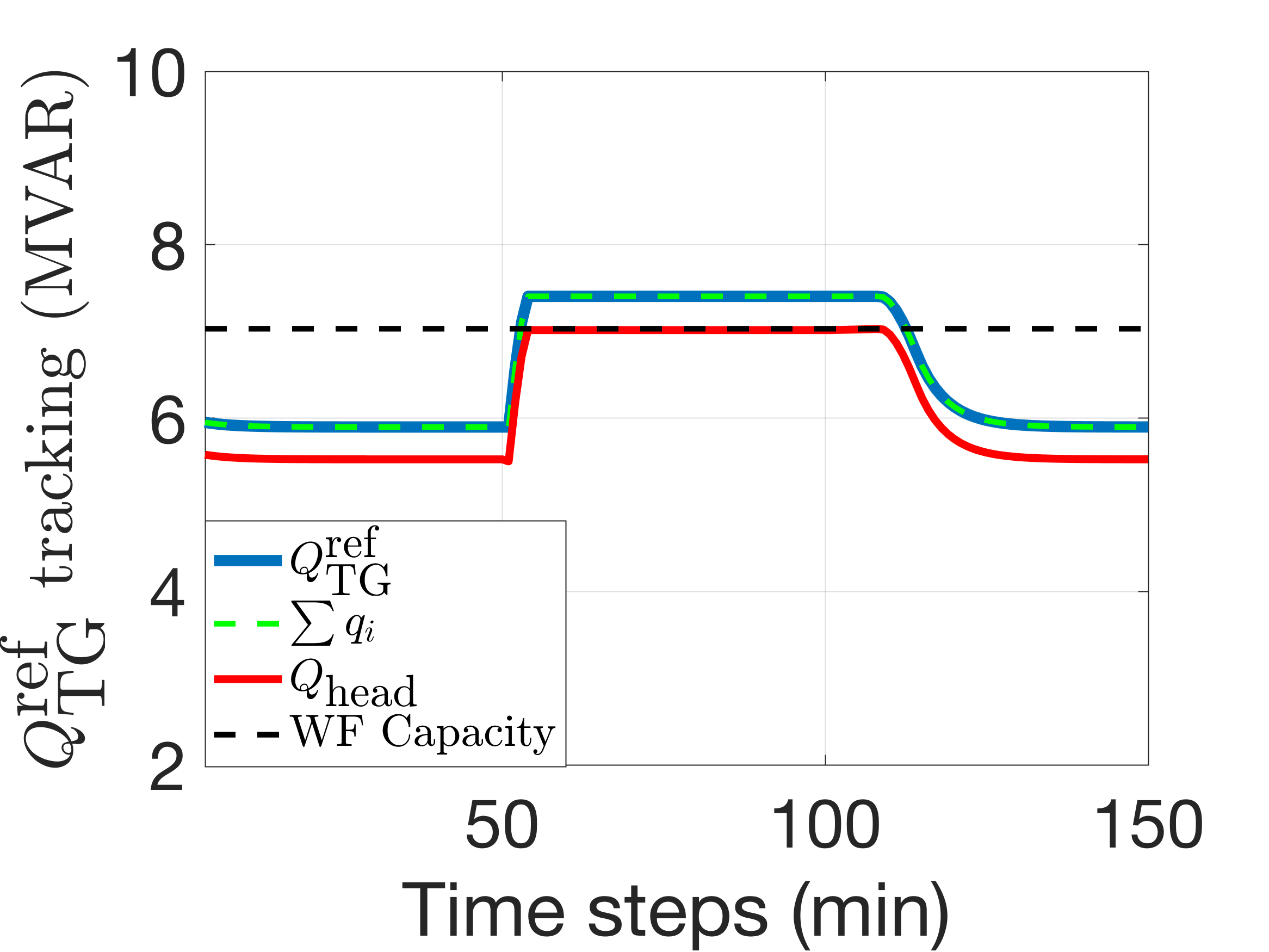}}
\hfill
\subfloat[\label{fig:Vwf_CIA}]{\includegraphics[width=0.33\linewidth]{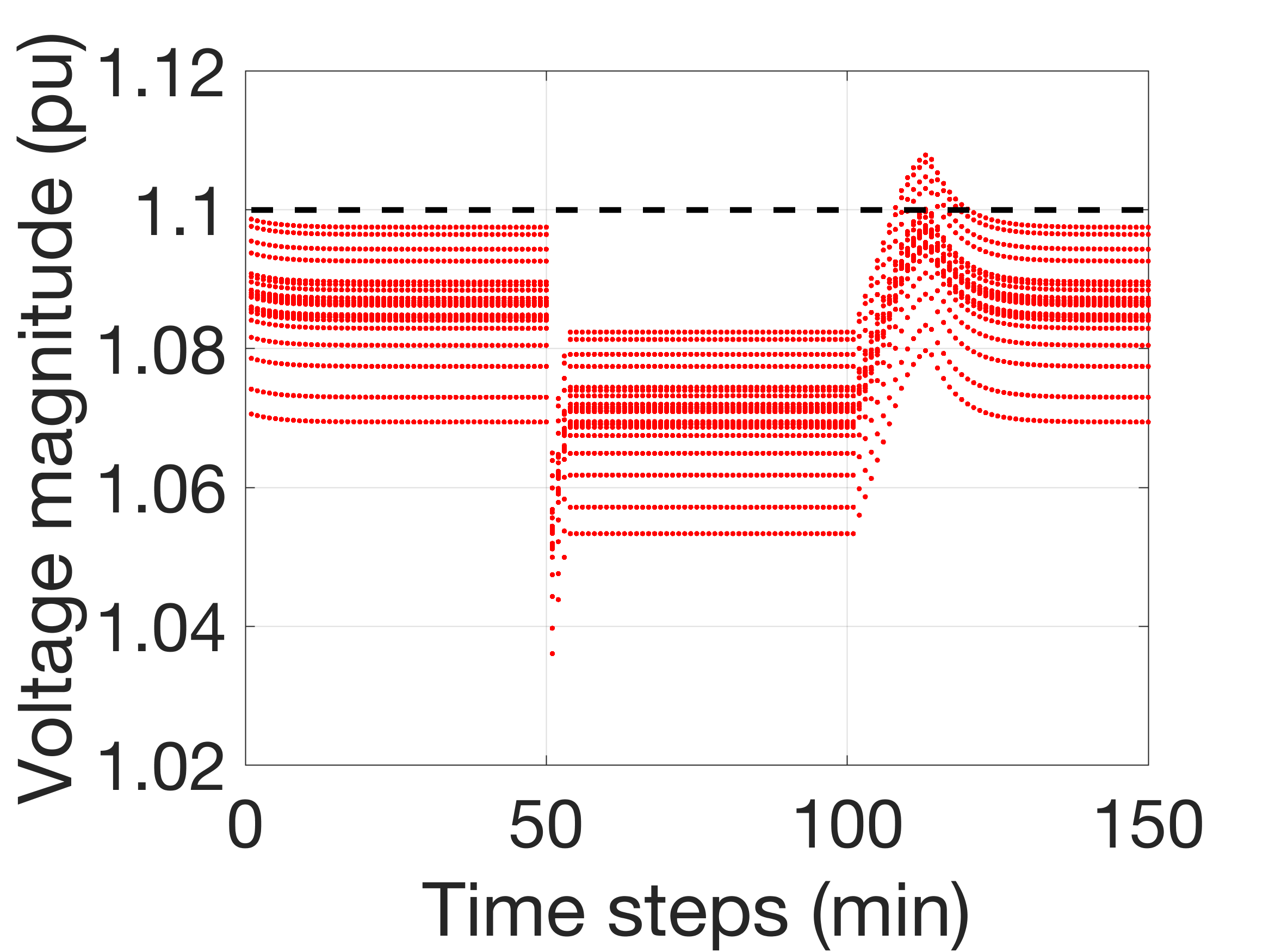}}
\\

% \subfloat
% [\label{fig:Qref_Martins}]{\includegraphics[width=0.33\linewidth]{images/Q_ref_tracking_Martins.png}}
% \hfill

% \subfloat
% [\label{fig:Vwf_Martins}]{\includegraphics[width=0.33\linewidth]{images/V_windfarm_Martins.png}}
% \hfill

\caption{Comparing the methods of disaggregating a desired reactive power set-point among the 19~wind turbines. (a,d): $V_{\text{PCC}}$ tracking under disturbance from grid voltage $V_{\text{grid}}$ for (a) grid-agnostic scheme (RMSE$=0.0055$ pu), and (d) the proposed CIA scheme (RMSE$=0.0113$ pu). (b,e): Head node physical reactive power $Q_{\text{head}}$ and the reference control signal for (b) grid-agnostic scheme, and (e) the proposed CIA scheme.  (c,f): wind farm network nodal voltages under (c) grid-agnostic scheme and (f) the proposed CIA scheme.}\label{fig:control_sim}
\end{figure*}

%Here we present a comparison of the developed control scheme with a scheme based on a grid-agnostic approach (does not consider the wind farm network constraints) and the scheme utilized in~\cite{Martins_2015}.
The different disaggregation schemes are compared in Fig.~\ref{fig:control_sim}. It is assumed that all the wind turbines are operating at their rated active power, i.e., we utilize scenario~1. Hence, for the method based on CIA, the nodal reactive power capacity is provided in Fig.~\ref{fig:scenario_1}. For the grid-agnostic approach, the network is ignored so the only constraints are the turbine reactive power limits ([-0.5,0.5]~MVAr in this case). Hence, for this method, the nodal reactive power capacity is the same as depicted in Fig.~\ref{fig:scenario_3}. The two methods for allocating reactive power across the wind turbines (grid-agnostic and CIA-based) result in the PCC voltages $V_{\text{PCC}}$ shown in Figs.~\ref{fig:Vref_grid_agnostic} and~\ref{fig:Vref_CIA}. Since the grid-agnostic approach employs the full range of nodal reactive capacities, the RMSE error in tracking $V_{\text{ref}}$ for the CIA-based approach is $\approx 2$ times that of the grid-agnostic approach. A comparison of the reactive power dispatches $Q_{\text{TG}}^{\text{ref}}$ (control signal) and $Q_{\text{head}}$ (physical quantity) is shown in Figs.~\ref{fig:Qref_grid_agnostic} and~\ref{fig:Qref_CIA}. From the figures it can be seen that the grid-agnostic approach provides much larger reactive power which results in the better voltage tracking performance. The CIA-based method provides reactive power up to the pre-calculated wind farm capacity (labelled ``WF Capacity'' in Fig.~\ref{fig:Qref_CIA}), which is obtained from the solution of (P1) and (P2). This results in the moderate voltage tracking performance. The difference between the dispatch control signal $Q_{\text{TG}}^{\text{ref}}$ and $Q_{\text{head}}$ is due to the losses in the network. The impact of the reactive power dispatch on the wind farm's  nodal voltages is depicted in Figs.~\ref{fig:Vwf_grid_agnostic} and~\ref{fig:Vwf_CIA}. It can be seen from Fig.~\ref{fig:Vwf_grid_agnostic} that the large reactive power output results in voltage violations across the wind farm. This is a consequence of the grid-agnostic approach not taking into account the wind farm network when dispatching wind turbine reactive power. As a result, it overestimates the network's reactive power capacity, resulting in voltage violations. In contrast, Fig.~\ref{fig:Vwf_CIA} shows that wind farm voltages are within limits, except for violations during transients. This is not unexpected as the CIA-based method only guarantees steady-state operating conditions. Extending the CIA-based method to capture transient behaviour is an interesting avenue for future work. This example illustrates the usefulness of the CIA-based method for enabling wind farms to provide reactive power support to transmission systems while ensuring reliability of the wind farm collector network.

% \begin{figure}[t]
% \includegraphics[width=0.4\textwidth]{images/V_ref_tracking_head.png}
% \caption{\label{fig:Vref_tracking_head}PCC voltage tracking response from the real-time control algorithm}
% \end{figure}

% \begin{figure}[t]
% \includegraphics[width=0.4\textwidth]{images/Q_ref_tracking_random.png}
% \caption{\label{fig:Qref_tracking}Tracking of the generated head-node reference reactive power signal along with the admissible range of operation.}
% \end{figure}

% \begin{figure}[t]
% \includegraphics[width=0.4\textwidth]{images/V_ref_tracking_random.png}
% \caption{\label{fig:Vref_tracking}Wind-farm voltages resulting from the real-time control algorithm}
% \end{figure}

\section{Conclusions and Future Work}\label{sec:concl}

The paper has considered the application of a convex inner approximation method for determining the reactive power support that can be provided by a wind farm. This method determines the nodal reactive power capacities that guarantee satisfaction of network constraints. These nodal capacities form the basis for a feedback control algorithm that provides voltage support to the grid by dispatching the reactive power of wind turbines. Through simulation of wind farm networks, we have shown the effectiveness of this approach.

Future work will consider systematic design of the objective function in (P1) and (P2), e.g., maximization of $Q_\text{head}$. We also seek to extend this work to incorporate system dynamic behaviour into the formulation of the CIA to guarantee transient operation of networks. Future work will also extend the method in~\cite{Martins_2015} to address lower reactive power limits and integrate that feedback methodology with the CIA-based method to improve performance for practical wind farms.
\textcolor{black}{Finally, we will seek to collaborate with wind industry experts to validate performance on high-fidelity models and support practical applications.}

\bibliographystyle{IEEEtran}
\small\bibliography{sample.bib}

\begin{thebibliography}{10}
\providecommand{\url}[1]{#1}
\csname url@samestyle\endcsname
\providecommand{\newblock}{\relax}
\providecommand{\bibinfo}[2]{#2}
\providecommand{\BIBentrySTDinterwordspacing}{\spaceskip=0pt\relax}
\providecommand{\BIBentryALTinterwordstretchfactor}{4}
\providecommand{\BIBentryALTinterwordspacing}{\spaceskip=\fontdimen2\font plus
\BIBentryALTinterwordstretchfactor\fontdimen3\font minus
  \fontdimen4\font\relax}
\providecommand{\BIBforeignlanguage}[2]{{%
\expandafter\ifx\csname l@#1\endcsname\relax
\typeout{** WARNING: IEEEtran.bst: No hyphenation pattern has been}%
\typeout{** loaded for the language `#1'. Using the pattern for}%
\typeout{** the default language instead.}%
\else
\language=\csname l@#1\endcsname
\fi
#2}}
\providecommand{\BIBdecl}{\relax}
\BIBdecl

\bibitem{slootweg2005wind}
J.~Slootweg, S.~De~Haan, H.~Polinder, and W.~Kling, ``Wind power and voltage
  control,'' \emph{Wind power in power systems}, pp. 413--432, 2005.

\bibitem{Camm_2009}
E.~H. {Camm}, M.~R. {Behnke}, O.~{Bolado}, M.~{Bollen}, M.~{Bradt},
  C.~{Brooks}, W.~{Dilling}, M.~{Edds}, W.~J. {Hejdak}, D.~{Houseman},
  S.~{Klein}, F.~{Li}, J.~{Li}, P.~{Maibach}, T.~{Nicolai}, J.~{Patino}, S.~V.
  {Pasupulati}, N.~{Samaan}, S.~{Saylors}, T.~{Siebert}, T.~{Smith},
  M.~{Starke}, and R.~{Walling}, ``Reactive power compensation for wind power
  plants,'' in \emph{2009 IEEE Power Energy Society General Meeting}, 2009, pp.
  1--7.

\bibitem{opila2010wind}
D.~F. Opila, A.~M. Zeynu, and I.~A. Hiskens, ``Wind farm reactive support and
  voltage control,'' in \emph{2010 IREP Symposium Bulk Power System Dynamics
  and Control-VIII (IREP)}.\hskip 1em plus 0.5em minus 0.4em\relax IEEE, 2010,
  pp. 1--10.

\bibitem{hiskens2013strategies}
I.~A. Hiskens, ``Strategies for voltage control and transient stability
  assessment,'' Univ. of Michigan, Ann Arbor, MI (United States), Tech. Rep.,
  2013.

\bibitem{Martins_2015}
J.~A. {Martin} and I.~A. {Hiskens}, ``Reactive power limitation due to
  wind-farm collector networks,'' in \emph{2015 IEEE Eindhoven PowerTech},
  2015, pp. 1--6.

\bibitem{silva2019loading}
V.~R.~N. Silva and R.~Kuiava, ``Loading margin sensitivity in relation to the
  wind farm generation power factor for voltage preventive control,''
  \emph{Journal of Control, Automation and Electrical Systems}, vol.~30, no.~6,
  pp. 1041--1050, 2019.

\bibitem{molzahn2017computing}
D.~K. Molzahn, ``Computing the feasible spaces of optimal power flow
  problems,'' \emph{IEEE Transactions on Power Systems}, vol.~32, no.~6, pp.
  4752--4763, 2017.

\bibitem{molzahn2019survey}
D.~K. Molzahn and I.~A. Hiskens, ``A survey of relaxations and approximations
  of the power flow equations,'' \emph{Now Publishers}, 2019.

\bibitem{nazir2019convex}
N.~Nazir and M.~Almassalkhi, ``Convex inner approximation of the feeder hosting
  capacity limits on dispatchable demand,'' in \emph{2019 IEEE 58th Conference
  on Decision and Control (CDC)}.\hskip 1em plus 0.5em minus 0.4em\relax IEEE,
  2019, pp. 4858--4864.

\bibitem{gan2014exact}
L.~Gan, N.~Li, U.~Topcu, and S.~H. Low, ``Exact convex relaxation of optimal
  power flow in radial networks,'' \emph{IEEE Transactions on Automatic
  Control}, vol.~60, no.~1, pp. 72--87, 2014.

\bibitem{brahma2020vb}
S.~Brahma, N.~Nazir, H.~Ossareh, and M.~R. Almassalkhi, ``Optimal and resilient
  coordination of virtual batteries in distribution feeders,'' \emph{IEEE
  Transactions on Power Systems}, vol.~36, no.~4, pp. 2841--2854, 2021.

\bibitem{nazir2019voltage}
N.~Nazir and M.~Almassalkhi, ``Voltage positioning using co-optimization of
  controllable grid assets in radial networks,'' \emph{IEEE Transactions on
  Power Systems}, vol.~36, no.~4, pp. 2761--2770, 2021.

\bibitem{nazir2019grid}
------, ``Grid-aware aggregation and realtime disaggregation of distributed
  energy resources in radial networks,'' \emph{IEEE Transactions on Power
  Systems}, pp. 1--1, 2021.

\bibitem{heidari2017non}
R.~Heidari, M.~M. Seron, and J.~H. Braslavsky, ``Non-local approximation of
  power flow equations with guaranteed error bounds,'' in \emph{Control
  Conference (ANZCC), 2017 Australian and New Zealand}, 2017, pp. 83--88.

\bibitem{zimmerman2011matpower}
R.~D. Zimmerman, C.~E. Murillo-S{\'a}nchez, and R.~J. Thomas, ``Matpower:
  Steady-state operations, planning, and analysis tools for power systems
  research and education,'' \emph{IEEE Transactions on power systems}, vol.~26,
  no.~1, pp. 12--19, 2011.

\end{thebibliography}

\end{document}